\documentclass[reqno]{amsart}
\usepackage{hyperref}
\usepackage{tikz}
\usepackage{amsmath}
\usepackage{epstopdf}
\usepackage{pdfpages}
  \usepackage{epsfig}

\begin{document}
\title[The multi-patch logistic equation with
asymmetric migration]{The multi-patch logistic equation with asymmetric migration}

\author[B.  ELBETCH, T. BENZEKRI, D.  Massart and T. Sari ]
{ELBETCH BILEL, BENZEKRI TOUNSIA, MASSART DANIEL, SARI TEWFIK}

\address{ELBETCH Bilel \newline
Department of Mathematics, University Dr. Moulay Tahar of Saida, Algeria.}
\email{bilel.elbetch@univ-saida.dz}

\address{BENZEKRI Tounsia \newline
Department of Mathematics, USTHB, Bab Ezzouar, Algiers, Algeria.}
\email{tbenzekri@usthb.dz}

\address{MASSART Daniel \newline 
IMAG, Univ Montpellier, CNRS, Montpellier, France.}
\email{daniel.massart@umontpellier.fr}

\address{SARI Tewfik, \newline 
ITAP, Univ Montpellier, INRAE, Montpellier SupAgro, Montpellier, France.}
\email{tewfik.sari@irstea.fr}

\subjclass[2010]{*****,*****,******}

\numberwithin{equation}{section}
\newtheorem{thm}{Theorem}
\numberwithin{thm}{section} 
\newtheorem{lm}[thm]{Lemma}
\newtheorem{defi}[thm]{Definition}
\newtheorem{con}[thm]{Conjecture}
\newtheorem{propo}[thm]{Proposition}
\newtheorem{cro}[thm]{Corollary}

\newtheorem{rem}[thm]{Remark}
\newtheorem{ex}[thm]{Example}
\newtheorem{tb}{Table}
\allowdisplaybreaks

\begin{abstract}
 This paper considers a multi-patch model, where each patch follows a logistic law, and patches are coupled by asymmetrical migration terms.
First, in the case of perfect
mixing,  i.e when the migration rate tends to infinity, the total population follows
a logistic law with a carrying capacity which in general is different from
the sum of the $n$ carrying capacities, and depends on the migration terms. Second, we determine,  in some particular cases, the conditions under which fragmentation and asymmetrical migration can lead to a total equilibrium population greater or smaller than the sum of the carrying capacities. Finally, for the three-patch model, we show numerically the existence of at least three critical values of the migration rate for which the total equilibrium population equals  the sum of the carrying capacities. 
\end{abstract}
\maketitle
\section{Introduction}
The study of the dynamics of a fragmented population is fundamental in theoretical ecology, with potentially very important applied aspects: 
what is the effect of migration on the general population dynamics ? What are the consequences of fragmentation on the persistence or extinction of the population ? 
When is a single large refuge  better or worse than several small ones (this is known as the SLOSS debate; see Hanski \cite{Hanski}) ?

The theoretical paradigm that has been used to treat these questions
is that of a single population fragmented into patches coupled by migration, and the sub-population in each patch follows a local logistic law. 
This system is modeled by a non linear system of differential equations of the following form:
\begin{equation}\label{4}
\frac{dx}{dt}=f(x)+\beta \Gamma x,
\end{equation}
where $x=(x_1,\ldots, x_n)^T$,  $n$ is the number of patches in the system, 
$x_{i}$ represents the population density in the $i$-th patch,
$f(x)=(f_{1}(x_{1}), \ldots, f_{n}(x_{n}))^T$, 
and 
\begin{equation}\label{logistics}
f_i(x_i)=r_ix_i(1-x_i/K_i),\quad i=1,\ldots n.
\end{equation} 
The parameters $r_{i}$ and $K_{i}$ are respectively the intrinsic growth rate and the carrying capacity of patch $i$.
The term $\beta \Gamma x$ on the right hand side of the system \eqref{4} describes the effect of the migration between the patches, where $\beta $ is the migration rate and 
$\Gamma=(\gamma_{ij})$ is the matrix representing the migrations  between the patches. For $i\neq j$, $\gamma_{ij}> 0$ denotes the incoming flux from patch $j$ to patch $i$.  
If $\gamma_{ij}=0$, there is no migration. The diagonal entries of $\Gamma$ satisfy the following equation
\begin{equation}\label{gammaiia}
\gamma_{ii}=-\sum_{j=1,j\neq i}^n\gamma_{ji},\qquad i=1,\cdots, n,
\end{equation}
which means that what comes out of a patch is distributed between the other  $n-1$ patches. 

In the absence of migration, ($\beta=0$), the system \eqref{4} admits $(K_{1},\ldots, K_{n})$  as a non trivial equilibrium point. 
This equilibrium is globally asymptotically stable (GAS) and the total population at equilibrium is equal to the sum of the carrying capacities. 
The problem is whether or not the equilibrium continues to be positive and GAS, for any $\beta>0$, and whether or not
the total population at equilibrium can be greater than the sum of the carrying capacities.
The case $n=2$ and $\Gamma$ symmetric
$$\Gamma=\left[
\begin{array}{rr}
-1&1\\1&-1
\end{array}
\right]$$
where $\gamma_{12}=\gamma_{21}$ is normalized to 1 has been considered by
Freedman and Waltman \cite{6} and  Holt \cite{9}. 
They analyzed the model in the case of perfect mixing $(\beta \rightarrow +\infty)$ and showed that the total equilibrium population can be greater than the sum of the carrying capacities $K_{1}+K_{2}$, so that patchiness has a beneficial effect on the total equilibrium population. More recently, Arditi et al. \cite{1} analyzed the behaviour of the system for all values of $\beta$. They  showed that only three situations occur: either  for any $\beta >0$, patchiness has a beneficial effect, or this effect is always detrimental, or the effect is beneficial for lower values of the migration coefficient $\beta$ and detrimental for  higher values. Arditi et al. \cite{2} extended these results to the case of two patches coupled by asymmetric migration, corresponding to the matrix
$$
\Gamma=\left[
\begin{array}{rr}
-\gamma_{21}&\gamma_{12}\\
\gamma_{21}&-\gamma_{12}
\end{array}
\right].
$$

DeAngelis et al. \cite{3,4} considered the case of $n > 2$ patches in a circle, with symmetric migration between any patch and its two neighbours : 
\begin{equation}\label{ModelDeAngelis}
\frac{dx_i}{dt}=r_ix_i\left(1-\frac{x_i}{K_i} \right)+\beta(x_{i-1}-2x_i+x_{i+1}),\qquad i=1,\ldots,n, 
\end{equation}
where we denote $x_0 = x_n$ and $x_{n+1} = x_1$, 
so that the same relationships hold between $x_i$, $x_{i-1}$ and $x_{i+1}$ for all values of $i$. 
This model corresponds to the matrix $\Gamma$ whose non-zero off-diagonal elements are given by
$$
\gamma_{1n}=\gamma_{n1}=1
\quad \mbox{ and } \quad
\gamma_{i,i-1}=\gamma_{i-1,i}=1, \quad \mbox{for} \quad 2\leq i\leq n.
$$
The system \eqref{ModelDeAngelis} is a one-dimensional discrete-patch version of the standard
reaction-diffusion model. In  \cite{3,4}  the perfect mixing case is described.

The case of the general symmetric migration was considered by the authors in \cite{4.1}. We studied  the system:
 \begin{equation}\label{01.3}
 \frac{d x_{i}}{dt}=r_ix_{i} \left( 1-\frac{x_i}{K_i}\right)+\beta \sum_{j=1, j\neq i}^{n}\gamma_{ij}(x_{j}-x_{i}), \hspace{1cm}i=1, \ldots, n,
 \end{equation}
where $\beta\gamma_{ij}$ is the rate of migration between patches $i$ and $j$. This system can be written in the form of System  \eqref{4} with 
$\Gamma=(\gamma_{ij})$, the symmetric matrix whose diagonal entries are defined  by \eqref{gammaiia}. 
We studied  the total population at equilibrium, as a function of the migration rate $\beta$. 
We gave conditions on the system parameters that ensure that  migration is beneficial or detrimental, and extended several results of \cite{1,3,4}. 

The aim of this work is to consider the case of  $n$ patches connected by asymmetric migration. Thus, we extend \cite{2} by considering the case $n\geq 2$, and we extend  \cite{4.1} by considering the case where $\Gamma$ is non symmetric. 

An important extension of \eqref{4} is the so called source-sink model, where the patches are of two types: 
the source patches, $1\leq i\leq m$, with logistic dynamics, and the sink patches, $m+1\leq i\leq n$, with exponential decay 
\begin{equation}\label{ED}
\left\{
\begin{array}{ll}
f_i(x_i)=r_ix_i(1-x_i/K_i),&i=1,\ldots,m,\\
f_i(x_i)=-r_ix_i,&i=m+1,\ldots,n. 
\end{array}
\right.
\end{equation}
The main problem is the number of source patches required for population persistence. 
For a recent study and bibliographical references the reader can consult Arino et al. \cite{arino} and Wu et al. \cite{wu}.
 
There is another important extension of (\ref{4},\ref{logistics}), where the dynamics on patch $i$ is of the form 
\begin{equation}\label{logisticGAO}
f_i(x_i)=r_ix_i(1-x_i/K_i)-\gamma_ix_i,\quad i=1,\ldots,n,
\end{equation}
with $\gamma_i>0$. This model is the limit system (when $t\to+\infty$) of an SIS model in $n$ patches connected by human migration. For details and further reading,  see Section \ref{SECSIS}. Note that, when  $r_i<\gamma_i$ for some patches, system (\ref{4},\ref{logisticGAO}) is  a source-sink model. Countrary to \eqref{ED}, the mortality in sink patch is density-dependent. For more details and bibliographical references the reader is referred to \cite{Gao}. 

Another example of source-sink model is the system considred by Nagahara et al. \cite{Nagahara}, called the ``island chain'' model, which is of the form:
\begin{equation}\label{Nagahara}
\frac{dx_i}{dt}=x_i\left(m_i-{x_i}\right)+\beta(x_{i-1}-2x_i+x_{i+1}),\qquad i=1,\ldots,n, 
\end{equation}
where we denote $x_0 = x_1$ and $x_{n+1} = x_n$.  This model is of the form \eqref{4},  $\Gamma$ being  the matrix  which verifies \ref{gammaiia}, and  whose non-zero off-diagonal elements are given by
$$
\gamma_{i,i-1}=\gamma_{i-1,i}=1, \quad \mbox{for} \quad 2\leq i\leq n.
$$
In the model \eqref{Nagahara} the ratios $\alpha_i=r_i/K_i$ in \eqref{logistics} are equal and are normalized to 1. The 
constant $m_i$ represents both the intrinsic
growth rate of the species in patch $i$ and the carrying capacity of the patch. If $m_i > 0$, then patch $i$ is favorable to the species. It is a source. The case $m_i=0$ is permitted and corresponds to a sink. The main purpose is to find the resource allocation   $(m_1,...,m_n)$ that maximizes the total population at equilibrium, under
the constraint that $\sum_im_i=m>0$ is fixed. For more details and information on the maximization of the total population with logistic growth in a patchy environment, the reader is referred to \cite{Nagahara} and the references therein.

For general information of the  effects of patchiness and migration in both continuous and discrete cases, and the results beyond the logistic model, the reader is referred to the work of Levin \cite{11,12}, DeAngelis et al. \cite{3,3.1,3.2,4},
Freedman et al. \cite{5}, Zaker et al. \cite{17}. 

It is worth noting that System \eqref{4} appears in metapopulation dynamics, involving explicit
movements of the individuals between distinct locations. For the graph theoretic and dynamical
system context in which metapopulation models are formulated, the reader is referred to Arino \cite[Section 2]{ArinoMetapopulation}.

The paper is organized as follows. In Section \ref{sec2}, the mathematical model of $n$ patches, and some preliminaries results, are introduced. 
In Section \ref{sec3}, the behavior of the model is studied when the migration rate tends to infinity. In Section \ref{sec4}, we compare the total equilibrium population with the sum of the carrying capacities in  some particular cases. In Section \ref{SECSIS}, the SIS patch model is considered, and the links with the logistic patch model are investigated. In Section \ref{sec5} the three-patch model is considered, and by numerical simulations we show the existence of a new behavior for the dynamics 
of the total equilibrium population as a function of the migration rate. 
In Appendix \ref{apa}, we recall some results for the two-patch model with asymmetrical migration. In Appendix \ref{apb}, we prove some useful auxiliary results.
\section{The mathematical model and preliminaries results}\label{sec2}
We consider the  model of multi-patch logistic growth, coupled by asymmetric migration terms 
\begin{equation}\label{m6}
\frac{dx_{i}}{dt}=r_{i}x_{i}\left(1-\frac{x_{i}}{K_{i}} \right) +
\beta\sum_{j=1, j\neq i} ^{n}\left(\gamma_{ij}x_{j}-\gamma_{ji}x_{i}\right), \qquad i=1,\cdots, n,
\end{equation}
where $\gamma_{ij}\geq 0$ denotes the
incoming flux from patch $j$ to patch $i$, for $i\neq j$. The system \eqref{m6} can be written in the form \eqref{4}, where $f$ is given by:  
\begin{equation}\label{f}
f(x)=\left(r_1x_1(1-x_1/K_1),\cdots,r_nx_n(1-x_n/K_n)\right)^T,
\end{equation}
and $\Gamma:=(\gamma_{ij})_{n\times n}$ is the matrix whose diagonal entries are given by \eqref{gammaiia}.
The matrix
$$
\Gamma_0:=\Gamma-{\rm diag}(\gamma_{11},\cdots,\gamma_{nn})
$$
which is the same as $\Gamma$, except that the diagonal elements are 0, 
is called the connectivity matrix. It is the adjacency matrix of the weighted directed graph $\mathcal{G}$, which has exactly $n$ vertices (the patches),  and has an arrow from patch $j$ to patch $i$, with weight $\gamma_{ij}$,  precisely when $\gamma_{ij} > 0$.

As to  the non-negativity of the solution, we have the following proposition:
\begin{propo}
The domain $\mathbb{R}_{+}^{n}=\left\lbrace (x_{1},\ldots, x_{n})\in \mathbb{R}^{n}/ x_{i}\geq 0, i=1,\ldots, n\right\rbrace $
is positively invariant for the system \eqref{m6}.
\end{propo}
\begin{proof}
The proof is the same as in the symmetrical case \cite[Prop 2.1]{4.1}.
\end{proof}

When the connectivity matrix $\Gamma_0$ is irreducible, System \eqref{m6} admits a unique  positive equilibrium $(x_{1}^{\ast}(\beta), \ldots, x_{n}^{\ast}(\beta))$, which is GAS,  see \cite[Theorem 2.2]{ArinoMetapopulation}, \cite[Theorem 1]{arino} or \cite[Theorem 6.1]{4.1}. 
In all of this work, we denote by 
$E^{\ast}(\beta)$ the positive equilibrium and by $X_{T}^{\ast}(\beta)$ the total population at equilibrium:
\begin{equation}\label{Ebeta}
E^{\ast}(\beta)=(x_{1}^{\ast}(\beta), \ldots, x_{n}^{\ast}(\beta)),
\qquad
X_{T}^{\ast}(\beta)=\sum_{i=1}^nx_{i}^{\ast}(\beta).
\end{equation}

\begin{rem}
The matrix $\Gamma_0$ being irreducible means that the weighted directed graph $\mathcal{G}$ is strongly connected, which means that every patch is reachable from every other patch, either directly or through other patches. The matrix $\Gamma$ is assumed to be irreducible throughout the rest of the paper.  
\end{rem}

\section{Perfect mixing}\label{sec3}
In this section our aim is to study the behavior of $E^{\ast}(\beta)$ and $X_{T}^{\ast}(\beta)$, defined by \eqref{Ebeta}, for large migration rate, i.e when $\beta\to\infty$. 
\subsection{The fast dispersal limit}
The following lemma was proved in  \cite[Lemma 2]{arino}, we include a proof for the ease of the reader.
\begin{lm}\label{lm41}
 Let $\Gamma$ be the migration matrix. Then, $0$ is a simple eigenvalue of $\Gamma$ and all non-zero eigenvalues of $\Gamma$ have negative real part. Moreover, 
the kernel  of the matrix $\Gamma$ 
is generated by a positive vector.

If the matrix $\Gamma$ is symmetric, then $\ker \Gamma$ is generated by $u=(1,...,1)^T $.
\end{lm}
\begin{proof}
Let $s=\max_{i=1,\ldots,n}(-\gamma_{ii} )$ and let $B$ be the matrix defined by
$$B=\Gamma+sI.$$
First, we note that since the matrix  $\Gamma$ verifies the property \eqref{gammaiia},  
then $\Gamma$  is a singular matrix and the vector  $u=(1,...,1)^T $ is an eigenvector of  $\Gamma ^T$ associated to the eigenvalue $0$. 
Thus  $u$  is an eigenvector of $B^T$, with eigenvalue $s$. 

The matrix $B^T$ is  non negative and irreducible, so by the Perron-Frobenius theorem  the spectral radius 
$$
\rho(B^T)=\max\left\{|\lambda|:\lambda \mbox{ is an eigenvalue of } B^T \right\},
$$
 is a simple eigenvalue of the matrix $B^T$ and it is the only eigenvalue of $B^T$ which admits a positive eigenvector, so $s= \rho(B^T)= \rho(B)$. 
 Therefore,  $\Gamma =B- \rho(B)I $ and $\dim(\ker \Gamma)=\dim(\ker \Gamma^T)=1$.
 
 All other eigenvalues of $B$ have modulus $< \rho (B)$, so their real parts are $< \rho (B)$. Since each eigenvalue of $\Gamma$ is $\lambda -  \rho (B)$, for some eigenvalue 
 $\lambda$ of $B$, all eigenvalues of $\Gamma$ have negative real part. 
 
Furthermore, according to the Perron-Frobenius theorem, there exists a positive vector $\delta $  such that $B \delta= \rho(B)\delta$, that is,  
 $\Gamma \delta=(B- \rho(B)I) \delta =0$. In particular, if the matrix  $\Gamma$   is symmetric then we may take $\delta = u$, that is, $\delta_i=1$, for all i.
\end{proof}
In all of this paper, we denote by $\delta=(\delta_1,\ldots, \delta_n)^T$ a positive vector which generates the vector space $\ker \Gamma$.

\begin{rem}\label{rem2}
The existence, uniqueness (mod. multiplicative factor), and positivity of $\delta$ were also proved  in Lemma 1 of Cosner et al. \cite{2.13}. 
On the other hand, it is shown in Guo et al. \cite[Lemma 2.1]{6.1} and Gao and Dong \cite[Lemma 3.1]{GaoDong} that the vector $(\Gamma_{11}^*,\ldots,\Gamma_{nn}^*)^T$ is a right eigenvector of $\Gamma$ associated with the
zero eigenvalue. Here, $\Gamma_{ii}^*$ is the cofactor of the $i$-th diagonal entry of  $\Gamma$. Therefore, we have   explicite formulae for the components of the vector $\delta$, as functions of  the
coefficients of $\Gamma$, at our disposal. 
For two patches we have $\delta=(\gamma_{12}, \gamma_{21})^T$, and for three patches we have $\delta=(\delta_1,\delta_2,\delta_3)^T$, where
\begin{equation}\label{coedel}
\left\lbrace 
\begin{array}{l}
\delta_1=\gamma_{12}\gamma_{13}+\gamma_{12}\gamma_{23}+\gamma_{32}\gamma_{13},\\
\delta_2=\gamma_{21}\gamma_{13}+\gamma_{21}\gamma_{23}+\gamma_{31}\gamma_{23},\\
\delta_3=\gamma_{21}\gamma_{32}+\gamma_{31}\gamma_{12}+\gamma_{31}\gamma_{32}.
\end{array}\right. 
\end{equation}

\end{rem}
The following result  asserts that when $\beta \rightarrow \infty$, the equilibrium $E^\ast(\beta)$ converges to an element of  $\ker \Gamma$.
\begin{thm}\label{pp2}
For the system \eqref{m6},  we have
$$
\lim_{\beta\to+\infty}E^{\ast}(\beta)=\dfrac{\sum_{i=1}^{n}\delta_{i}r_{i}}{\sum_{i=1}^{n}\delta_{i}^2\alpha_{i}}\left(\delta_{1}, \ldots, \delta_{n}  \right),
$$
where $\alpha_i=r_i/K_i$.
\end{thm}

\begin{proof}
Denote
$$E^{\ast}(\infty)=\left(\delta_{1}\dfrac{\sum_{i=1}^{n}\delta_{i}r_{i}}{\sum_{i=1}^{n}\delta_{i}^2\alpha_{i}}, \ldots, \delta_{n} \dfrac{\sum_{i=1}^{n}\delta_{i}r_{i}}{\sum_{i=1}^{n}\delta_{i}^2\alpha_{i}}  \right).$$	
Dividing Equation \ref{4} at the equilibrium $E^{\ast}(\beta)$ by $\beta$, for $\beta >0$, yields
\[ 
\mbox{for all } \beta >0, \ \frac{1}{\beta} f(E^{\ast}(\beta))+ \Gamma E^{\ast}(\beta)=0.
 \]	
Thus any limit point, when $\beta \rightarrow \infty$,  of the set $\left\lbrace E^{\ast}(\beta) : \beta >0 \right\rbrace $ lies in the kernel of $\Gamma$. 
Now, taking the sum of all equations in
\[
r_{i}x_{i}\left(1-\dfrac{x_{i}}{K_{i}} \right) +
\beta\sum_{j=1, j\neq i} ^{n}(\gamma_{ij}x_{j}-\gamma_{ji}x_{i})=0, 
\qquad 
i=1,\cdots, n,
\]

 we see that $E^{\ast}(\beta)$ lies in the ellipsoid 
$$
\mathbb{E}^{n-1}=\left\{x\in\mathbb{R}^n:
\Theta(x):=\sum_{i=1}^{n}r_{i}x_{i}\left( 1-\dfrac{x_{i}}{K_{i}}\right)=0\right\}.
$$
The ellipsoid $\mathbb{E}^{n-1}$ is compact, so the equilibrium $E^\ast(\beta)$ has at least one limit point in $\mathbb{E}^{n-1}$, when $\beta$ goes to infinity. Since the kernel of $\Gamma$ has dimension $1$, and 
$\mathbb{E}^{n-1}$ is the boundary of a convex set, $\mathbb{E}^{n-1} \cap \ker \Gamma$ consists of at most two points. Since the origin  and $E^\ast(\infty)$ both lie in 
$\mathbb{E}^{n-1} \cap \ker \Gamma$, we get that 
\[
\mathbb{E}^{n-1} \cap \ker \Gamma = \left\lbrace 0, E^\ast(\infty)\right\rbrace .
\]

Therefore, to prove the convergence of $E^\ast(\beta)$ to
$E^\ast(\infty)$, it suffices to prove that the origin cannot be a limit point of $E^\ast(\beta)$. We claim that for any $\beta$, there exists $i$ such that $x_i^\ast(\beta)\geq K_i$, which entails that $E^\ast(\beta)$ is bounded away from the origin. The coordinates of the vector $\Gamma E^\ast(\beta)$ sum to zero, hence at least one of them, say, the $i$-th, is non negative. Then
$$r_ix_i^\ast(\beta)\left(1-\dfrac{x_i^\ast(\beta)}{K_i}\right) \leq 0,$$
and since $x_i^\ast(\beta)$ cannot be negative or $0$, we have $x_i^\ast(\beta)\geq K_i$.

\end{proof}
As a corollary of the previous theorem, we obtain the following result, which describes the total equilibrium population for perfect mixing:

\begin{propo}\label{p4.2}
We have 
\begin{equation}\label{5151}
X_{T}^{\ast}(+\infty)=\lim_{\beta\to+\infty}\sum_{i=1}^{n}x_{i}^{\ast}(\beta)=\left(\sum_{i=1}^{n}\delta_i\right)   \dfrac{\sum_{i=1}^{n}\delta_i r_{i}}{\sum_{i=1}^{n}\delta_i^{2}\alpha_{i}}.
\end{equation}
Denote $K=(K_1,\ldots, K_n)^T$. If $K=\lambda \delta$ with $\lambda>0$, that is to say $K\in \ker \Gamma$, then
$
X_{T}^{\ast}(+\infty)=\lambda \sum_{i=1}^{n}\delta_i=\sum_{i=1}^{n}K_i.
$
\end{propo}
\begin{proof}
For the proof of \eqref{5151}, it suffices to sum  the $n$ components of the point $E^{\ast}(\infty)$. For the  case $K\in\ker\Gamma$, it suffices to
replace $K_i$ by $\lambda \delta_i$ in \eqref{5151}.
\end{proof}

Actually, when $K\in \ker\Gamma$, we have $X_T^\ast(\beta)=\sum_i K_i$ for all $\beta>0$, see Prop. \ref{p52}.

In the case $n=2$, one has $\delta_1=\gamma_{12}$ and $\delta_2=\gamma_{21}$, as shown in Remark \ref{rem2}. Therefore \eqref{5151} becomes
$$
X_T^* (+\infty) = (\gamma_{12}+\gamma_{21})
\frac{\gamma_{12}r_1+\gamma_{21}r_2}
{\gamma_{12}^2\alpha_1+\gamma_{21}^2\alpha_2},
$$
which is the formula \cite[Equation 7]{2} given by Arditi et al. 

If the matrix $\Gamma$ is symmetric, one has $\delta_i=1$, for all $i$, as shown in Lemma \ref{lm41}. 
Therefore \eqref{5151}  specializes to the formula given in \cite[Equation (24)]{4.1}:
$$
X_T^\ast(+\infty)=n\dfrac{\sum_{i=1}^n r_i}{\sum_{i=1}^n r_i/K_i}.
$$

\subsection{Two time scale dynamics}
In  \cite{4.1} the authors also obtained the formula \eqref{5151}, in the symmetrical n-patch  case (i.e the matrix $\Gamma$ is symmetric), by using singular perturbation theory, see \cite[Theorem 4.6]{4.1}. 

We showed that, if $\left(x_1(t,\beta),\ldots, x_n(t,\beta)\right)$ is the solution of
\eqref{01.3}, with initial condition $(x_{1}^0,\ldots,x_{n}^0)$, 
then, when $\beta\to\infty$, 
the total population $\sum x_i(t,\beta)$ is approximated by $X(t)$, the solution of the  logistic equation
\begin{equation}\label{reduit1}
\frac{dX}{dt}=rX\left(1-\frac{X}{nK}\right),
\mbox{ where }r=\frac{\sum _{i=1}^n r_i}{n},\quad K=\frac{\sum _{i=1}^n r_i}{\sum _{i=1}^n\alpha_i}
\mbox{ and }\alpha_i=\frac{r_i}{K_i}
\end{equation}
with initial condition $X_0=\sum x_{i}^0$. Therefore the total population behaves like the solution of the  logistic equation given by \eqref{reduit1}.
In addition, one obtains the following property: 
with the exception of a small initial interval, the population densities  $x_i(t,\beta)$ are  approximated by $X(t)/n$, see \cite[Formula (37)]{4.1}. 
Therefore, this approximation shows that, when $t$ and $\beta$ tend to $\infty$, the population density  $x_i(t,\beta)$ tends toward $\frac{\sum r_i}{\sum \alpha_i}$, and in addition,  $x_i(t,\beta)$ quickly jumps from its initial condition $x_{i}^0$ to the average $X_0/n$ and then is very close to 
$X(t)/n$.   
Our aim is to generalize this result  for the asymmetrical $n$-patch model \eqref{m6} (i.e the matrix $\Gamma$ is non symmetric). We have the following result

\begin{thm}\label{th71}
Let $(x_{1}(t,\beta),\ldots, x_{n}(t,\beta))$ be the solution of the system \eqref{m6} with initial condition $(x_{1}^{0},\cdots, x_{n}^{0})$ satisfying $x_{i}^{0}\geq 0$ for $ i=1\cdots n$.
Let $Y(t)$ be the solution of the logistic equation 
\begin{equation}\label{reduitneq}
\frac{dX}{dt}=rX\left(1-\dfrac{X}{\left[\sum_{i=1}^{n}\delta_i \right]  K}\right),
\end{equation}
where
\begin{equation}\label{reduitn}
r=\dfrac{\sum_{i=1}^n \delta_i r_i}{\sum_{i=1}^{n}\delta_i }, 
K=\dfrac{\sum_{i=1}^n \delta_i r_i}{\sum_{i=1}^n \delta_i^{2} \alpha_i}
\mbox{ and }\alpha_i=\frac{r_i}{K_i},
\end{equation}
with initial condition $X_0=\sum_{i=1}^n x_{i}^0$. Then,
when $\beta\to\infty$, we have 
\begin{equation}\label{ax}
\sum_{i=1}^nx_{i}(t,\beta)={Y(t)}+o(1),\qquad  \mbox{ uniformly for } t\in[0,+\infty)
\end{equation}
and, for any $t_0>0$,  we have
\begin{equation}\label{ax_i}
x_{i}(t,\beta)=\frac{\delta_i }{\sum_{i=1}^{n}\delta_i }Y(t)+o(1), \quad i=1,\ldots, n,\mbox{ uniformly for }\quad  t\in[t_0,+\infty).
\end{equation}
\end{thm}
\begin{proof}
Let $X(t,\beta)=\sum_{i=1}^n x_{i}(t,\beta)$. We rewrite the system \eqref{m6} using the variables $(X, x_{1}, \cdots, x_{n-1})$, and get:
\begin{equation}\label{43}
\left\lbrace 
\begin{array}{rcl}
\dfrac{dX}{dt}&=&
\displaystyle\sum_{i=1}^nr_{i}x_{i}\left(1-\dfrac{x_{i}}{K_{i}} \right),\\
\dfrac{dx_{i}}{dt}&=&r_{i}x_{i}\left(1-\dfrac{x_{i}}{K_{i}} \right) +
\beta\displaystyle\sum_{j=1, j\neq i} ^{n}(\gamma_{ij}x_{j}-\gamma_{ji}x_{i}), 
\quad i=1,\cdots, n-1.\\
\end{array}\right. 
\end{equation}
This system is actually a system in the variables $(X, x_{1}, \cdots, x_{n-1})$, since, whenever $x_n$ appears in the right hand side of \eqref{43}, it should be replaced by
\begin{equation}\label{xn}
x_n=X-\sum_{i=1}^{n-1}x_i.
\end{equation}
When $\beta\to\infty$, \eqref{43} is a \textit{slow-fast} system, with one \textit{slow variable}, $X$, and $n-1$ \textit{fast variables}, $x_{i}$ for $i=1\cdots n-1$.
As suggested by Tikhonov's theorem \cite{10,15,16}, we consider the dynamics of the fast variables in the time scale $\tau=\beta t$. We get
$$
\dfrac{dx_{i}}{d\tau}=\dfrac{1}{\beta}r_{i}x_{i}\left( 1-\dfrac{x_{i}}{K_{i}}\right)+\displaystyle\sum_{j=1, j\neq i} ^{n}(\gamma_{ij}x_{j}-\gamma_{ji}x_{i}), 
\quad i=1,\cdots, n-1. 
$$
where $x_n$ is given by \eqref{xn}.
In the limit $\beta \to \infty$, we find the \textit{fast dynamics}
$$
\dfrac{dx_{i}}{d\tau}=\displaystyle\sum_{j=1, j\neq i} ^{n}(\gamma_{ij}x_{j}-\gamma_{ji}x_{i}), 
\qquad i=1,\cdots, n-1.
$$
This is an ($n-1$)-dimensional linear differential system in the variable $Z:=(x_1,\cdots,x_{n-1})$, which can be rewritten in matricial form:
\begin{equation}\label{0410}
\dot{Z}=\mathcal{L}Z+XV, \quad \text{ with }\qquad \mathcal{L}:=L-U,
\end{equation}
 where 
 $ L:=(\gamma_{ij})_{n-1 \times n-1}$ is the sub matrix of the matrix $\Gamma$, 
 obtained by dropping the last row and the last column of $\Gamma$, $V$ is the vector defined by $V:=(\gamma_{in})_{n-1 \times 1}$ and $U=(V;\ldots; V)$.

By Lemma \ref{appBlm1},  the matrix $\mathcal{L}$ is stable, that is, all of its eigenvalues have negative real part.
Therefore, it is invertible and the equilibrium of the system \eqref{0410}  is GAS.
This equilibrium is given by
$$
 \left( \dfrac{\delta_{1}}{\sum_{i=1}^{n}\delta_{i}}X, \ldots, \dfrac{\delta_{n-1}}{\sum_{i=1}^{n}\delta_{i}}X\right) ^{T}.
 $$
 Indeed, we denote by $L^{(i)}, U^{(i)}$ and $V^{(i)}$ the i-th row of the matrix $L, U$ and the vector $V$ respectively. We have:
\begin{align*}
\dfrac{\delta_{n}}{\sum_{i=1}^{n}\delta_{i}}\left( L^{(i)}-U^{(i)}\right) \left( 
\begin{array}{ccc}
\dfrac{\delta_{1}}{\delta_{n}}X&\ldots&
\dfrac{\delta_{n-1}}{\delta_{n}}X
\end{array}\right)^{T}&=-\dfrac{\delta_{n}}{\sum_{i=1}^{n}\delta_{i}}X\gamma_{in}-\dfrac{\sum_{i=1}^{n-1}\delta_{i}}{\sum_{i=1}^{n}\delta_{i}}X\gamma_{in}\\
&=-X\gamma_{in}=-XV^{(i)}.
\end{align*}
%
Thus,  the slow manifold of System \eqref{43} is given by 
 \begin{equation}\label{slowm}
 x_i= \dfrac{\delta_{i}}{\sum_{i=1}^{n}\delta_{i}}X, \qquad i=1, \ldots, n-1.
 \end{equation}
 As this manifold is {GAS}, Tikhonov's theorem  ensures that after a fast transition toward the slow manifold, the solutions of  \eqref{43} are approximated by the solutions of the \textit{reduced model}, which is obtained by replacing \eqref{slowm} into the dynamics of the slow variable, that is:
$$
\frac{dX}{dt}=
\sum_{i=1}^n r_{i}\dfrac{X}{\sum_{i=1}^{n}\delta_i}\delta_i \left(1-\dfrac{X}{\left(\sum_{i=1}^{n}\delta_i\right)  K_{i}}\delta_i \right)=rX\left(1-\frac{X}{\left(\sum_{i=1}^{n}\delta_i\right)  K}\right),
$$
where $r$ and $K$ are defined in  \eqref{reduitn}. Therefore, the reduced model is \eqref{reduitneq}. 
Since \eqref{reduitneq} admits 
$$X^\ast=\left(\sum_{i=1}^{n}\delta_i\right)  K=\left(\sum_{i=1}^{n}\delta_i\right)\frac{\sum_{i=1}^{n} \delta_i r_i}{\sum_{i=1}^{n} \delta_i^2 \alpha_i}$$
as a positive equilibrium point, which is {GAS} in the positive axis, 
the approximation given by Tikhonov's theorem holds for all $t\geq 0$ for the slow variable and for all $t\geq t_0>0$ for the fast variables, where $t_0$ is as small as we want. Therefore, letting $Y(t)$ be the solution of the reduced model
\eqref{reduitneq} with initial condition $Y(0)=X(0,\beta)=\sum_{i=1}^n x_{i}^0$, then,  
when $\beta \to\infty$, we have the approximations \eqref{ax} and \eqref{ax_i}.
\end{proof}

In the case of perfect mixing, the approximation \eqref{ax} shows that the total population behaves like the solution of the single logistic equation \eqref{reduitn} and then,  
when $t$ and $\beta$ tend to $\infty$, the total population $\sum x_i(t,\beta)$ tends toward 
$\left(\sum_{i=1}^{n}\delta_i\right)K=\left(\sum_{i=1}^{n}\delta_i\right)\frac{\sum \delta_i r_i}{\sum \delta_i^2\alpha_i}$ as stated in Prop. \ref{p4.2}. 
The approximation \eqref{ax_i} shows that, with the exception of a thin initial boundary layer, where the population density  
$x_i(t,\beta)$ quickly jumps from its initial condition $x_{i}^0$ to $
\delta_i X_0/\sum_{i=1}^{n}\delta_i$,
each patch of the n-patch model behaves like the logistic equation
\begin{equation}\label{328}
\frac{du}{dt}=ru\left(1-\frac{u}{\delta_i K}\right)\mbox{ where}\quad r=\dfrac{\sum_{i=1}^n \delta_i r_i}{\sum_{i=1}^{n}\delta_i },\quad 
K=\dfrac{\sum_{i=1}^n \delta_i r_i}{\sum_{i=1}^n \delta_i^{2} \alpha_i},\quad
\alpha_i=\frac{r_i}{K_i}.
\end{equation} 
Hence, when $t$ and $\beta$ tend to $\infty$, the population density  $x_i(t,\beta)$ tends toward $\delta_i \frac{\sum \delta_i r_i}{\sum \delta_i^2\alpha_i}$, 
as stated in Theorem \ref{pp2}.
\begin{rem}
The single logistic equation \eqref{328} gives an approximation of the population density in each patch in the case of perfect mixing. 
The intrinsic growth rate $r$ in \eqref{328} is the  arithmetic mean of the $r_1,\ldots, r_n,$ weighted by $\delta_1,\ldots,\delta_n$, 
and the carrying capacity $K$ is the harmonic mean of  $K_i/\delta_i$, weighted by $\delta_i r_i, i=1,\ldots, n$. 
We point out the similarity between our expression for the carrying capacity in the limit $\beta \to \infty$, and the expression obtained in spatial homogenization, 
see e.g \cite[Formula 81]{Yurk} and also \cite[Formula 28]{17}.

\end{rem}
\subsection{Comparison of $X_T^*(+\infty)$ with $\sum_iK_i$. }
According to Formula \eqref{5151}, it is  clear that the total equilibrium population at $\beta=0$ and at $\beta=+\infty$ are different in general.

 In the remainder of this section, we give some conditions, in  the space of parameters $r_i, K_i, \alpha_i$ and $\delta_i$,  for limit of  the total equilibrium population when 
 $\beta\to \infty$ to be greater or smaller than the sum of the carrying capacities. We show that all three cases are possible, i.e $X_T^\ast(+\infty)$ can be greater than, smaller than, 
  or equal to $X_T^\ast(0)$. 
First, we start by giving some particular values  of the parameters for which  equality holds.
\begin{propo}\label{3.9}
Consider the system \eqref{m6}.
 If the vector $\left( \frac{1}{\alpha_1},\ldots, \frac{1}{\alpha_n}\right)^T$ lies in $\ker \Gamma$, then $X_T^\ast(+\infty)=\sum_iK_i$.
\end{propo}
\begin{proof}
Direct consequence of the equation \eqref{5151}.
\end{proof}
Note that, if the matrix $\Gamma$ is symmetric, then by Lemma \ref{lm41}, Prop. \ref{3.9} says that if all $\alpha_i$ are equal, then $X_T^\ast(\infty)=\sum_iK_i$, which is \cite[Prop 4.4]{4.1}.

In the next proposition, we give two cases which ensure that $X_T^\ast(0)$ can be greater or smaller than $X_T^\ast(+\infty)$. This result can be stated as the following proposition:
\begin{propo}\label{prop12}
Consider the system \eqref{m6}.
\begin{enumerate}
\item If $\dfrac{K_1}{\delta_1}\leq \ldots \leq \dfrac{K_n}{\delta_n} \mbox{ and } \delta_1 \alpha_1 \leq \ldots \leq \delta_n \alpha_n$, 
or if $\dfrac{K_1}{\delta_1}\geq \ldots \geq \dfrac{K_n}{\delta_n}$ and $\delta_1 \alpha_1 \geq \ldots \geq \delta_n \alpha_n$,
then $X_T^\ast(+\infty)\geq X_T^\ast(0)$.
\item If $\dfrac{K_1}{\delta_1}\geq \ldots \geq \dfrac{K_n}{\delta_n}$ and $\delta_1 \alpha_1 \leq \ldots \leq \delta_n \alpha_n$, 
or if $\dfrac{K_1}{\delta_1}\leq \ldots \leq \dfrac{K_n}{\delta_n}$ and $\delta_1 \alpha_1 \geq \ldots \geq \delta_n \alpha_n$,
then $X_T^\ast(+\infty)\leq X_T^\ast(0) $.
\end{enumerate}
In both items, if at least one of the inequalities in $\dfrac{K_1}{\delta_1}\leq \ldots \leq \dfrac{K_n}{\delta_n}$ or $\dfrac{K_1}{\delta_1}\geq \ldots \geq \dfrac{K_n}{\delta_n} $
is strict, then the inequality is strict in the conclusion. 
\end{propo}
\begin{proof}
Apply Lemma \ref{lma41} with  the following choice:
$w_i=\delta_i$, $u_i=\dfrac{K_i}{\delta_i}$, and $v_i=\delta_i \alpha_i$, for all $ i=1, \ldots, n$.
\end{proof}
If the matrix $\Gamma$ is symmetric, one has $\delta_i = 1$, for all  $i$, as shown in Lemma \ref{lm41}. Therefore Prop. \ref{prop12} becomes
 
\begin{cro}
 Consider the system \eqref{m6}. Assume that $\Gamma$ is symmetric.
\begin{enumerate}
\item If $K_1\leq \ldots \leq K_n \mbox{ and } \alpha_1 \leq \ldots \leq  \alpha_n$, 
or if $K_1\geq \ldots \geq K_n$ and $ \alpha_1 \geq \ldots \geq  \alpha_n$,
then $X_T^\ast(+\infty)\geq X_T^\ast(0)$.
\item If $K_1\geq \ldots \geq K_n$ and $ \alpha_1 \leq \ldots \leq \alpha_n$, 
or if $K_1 \leq \ldots \leq K_n$ and $ \alpha_1 \geq \ldots \geq  \alpha_n$,
then $X_T^\ast(+\infty)\leq X_T^\ast(0) $.
\end{enumerate} 
\end{cro}
 
This result implies Items 1 and 2 of  \cite[Theorem B.1]{3.2}, which were obtained for the model  (\ref{ModelDeAngelis}) in the particular case $r_i = K_i$.



\section{Influence of asymmetric dispersal on total population size}\label{sec4}

In this section, we will compare, in some particular cases of the system \eqref{m6}, the total equilibrium population 
$X_{T}^{\ast}(\beta)=x_{1}^{\ast}(\beta)+ \ldots+x_{n}^{\ast}(\beta)$,
 with the sum of carrying capacities denoted by ${X}_{T}^\ast(0)=K_{1}+\ldots+ K_{n}$,
when the rate of migration $\beta$ varies from zero to infinity. We show that the total equilibrium population, 
$X_T^\ast(\beta)$, is generally different from the sum of the carrying capacities ${X}_{T}^\ast(0)$. 
Depending on the local parameters of the patches and the kernel of the matrix $\Gamma$, $X_T^\ast(\beta)$ 
can either be greater than, smaller than,  or equal to the sum of the carrying capacities.
\subsection{Asymmetric dispersal may be unfavorable to the total equilibrium population}
When $\Gamma$ is symmetric, we have already proved that if all the growth rates 
are equal then dispersal is always unfavorable to the total equilibrium 
population, see \cite[Prop. 3.1]{4.1}. We   also noticed that the result still 
holds in the general case when $\Gamma$ is not necessarily symmetric, see \cite[Prop. 6.2]{4.1}. 
Hence we have the following

\begin{propo}
 If $r_1=\ldots=r_n$ then
\[
X_{T}^{\ast}(\beta)=\sum_{i=1}^{n}x_{i}^{\ast}(\beta) \leq \sum_{i=1}^{n}K_{i},  \qquad\text{ for all } \beta \geq 0.
\]
\end{propo}

For a two-patch logistic model, this result has been proved by Arditi et al. \cite[Prop. 2, item 3]{1}
for symmetric dispersal  and for asymmetric 
dispersal \cite[Prop. 1, item 3]{2}.

\subsection{Asymmetric dispersal may be favorable to the total equilibrium population}
In this section, we give a situation where the dispersal is favorable to the total equilibrium population. Mathematically speaking:
\begin{propo}\label{p55}
 Assume that for all $j<i$, $ \alpha_{i}\gamma_{ij}=\alpha_{j}\gamma_{ji}$. Then
$$
X_{T}^{\ast}(\beta)\geq \sum_{i=1}^{n}K_{i}\qquad \text{ for all } \beta \geq 0.
$$
Moreover, if there exist $i_0$ and $j_0\neq i_0$ such that $r_{i_0}\neq r_{j_0} $, then $X_{T}^{\ast}(\beta)> \sum_{i=1}^{n}K_{i},$  for all $  \beta > 0.$
\end{propo}
\begin{proof}
The equilibrium point $E^{\ast}(\beta)$ satisfies the system
\begin{equation}\label{19}
0=\alpha_{i}x^{\ast}_{i}(\beta)\left( K_{i}-x^{\ast}_{i}(\beta)\right) +\beta \sum_{j=1 , j \neq i}^{n}(\gamma_{ij} x^{\ast}_{j}(\beta)- \gamma_{ji}x^{\ast}_{i}(\beta)), \quad i=1\cdots n.
\end{equation}
Dividing \eqref{19} by $\alpha_{i}x_{i}^{\ast}$, one obtains
$$
x^{\ast}_{i}(\beta)=
K_{i}+\beta \sum_{j=1 , j \neq i}^{n}\dfrac{\gamma_{ij}x_{j}^{\ast}(\beta)-\gamma_{ji}x_{i}^{\ast}(\beta)}{\alpha_{i}x_{i}^{\ast}(\beta)}
.$$
Taking the sum of these expressions shows that the total equilibrium population $X_{T}^{\ast}$ satisfies the following relation: 
\begin{align}\label{2424}
X_{T}^{\ast}(\beta)&=\sum_{i=1}^{n}K_{i}
+\beta \sum_{i=1}^{n}\sum_{j=1 , j \neq i}^{n}\dfrac{\gamma_{ij}x_{j}^{\ast}(\beta)-\gamma_{ji}x_{i}^{\ast}(\beta)}{\alpha_{i}x_{i}^{\ast}(\beta)}
\nonumber
\\ \nonumber
&=\sum_{i=1}^{n}K_{i}+\beta\sum_ {j<i}\left( \dfrac{\gamma_{ij}x_{j}^{\ast}(\beta)-\gamma_{ji}x_{i}^{\ast}(\beta)}{\alpha_{i}x_{i}^{\ast}(\beta)}+\dfrac{\gamma_{ji}x_{i}^{\ast}(\beta)-\gamma_{ij}x_{j}^{\ast}(\beta)}{\alpha_{j}x_{j}^{\ast}(\beta)}\right) 
\\ 
&=\sum_{i=1}^{n}K_{i}+\beta\sum_ {j<i}\dfrac{\left(\gamma_{ij}x_{j}^{\ast}(\beta)-\gamma_{ji}x_{i}^{\ast}(\beta)\right)\left(\alpha_{j}x_{j}^{\ast}(\beta)-\alpha_{i}x_{i}^{\ast}(\beta)\right)}{ \alpha_{j}\alpha_{i}x_{j}^{\ast}(\beta)x_{i}^{\ast}(\beta)}.
\end{align}
The conditions $\alpha_i\gamma_{ij}=\alpha_j\gamma_{ji}$ can be written  
$\kappa_{ij}:=\alpha_i/\gamma_{ji}=\alpha_j/\gamma_{ij}$ for all $j<i$, such that $\gamma_{ij}\neq 0$  and $\gamma_{ji}\neq 0$. Therefore, there exists $\kappa_{ij}>0$ such that 
$$
\alpha_{j}=\kappa_{ij}\gamma_{ij}
\text{ and }
\alpha_{i}=\kappa_{ij}\gamma_{ji} 
\text{ for all } i , j \text{ with } \gamma_{ij}\neq 0 \text{ and } \gamma_{ji}\neq 0.
$$
Replacing $\alpha_i$ and $\alpha_j$ in \eqref{2424}, one obtains
\begin{equation}\label{som}
X_{T}^{\ast}(\beta)=\sum_{i=1}^{n}K_{i}
+\beta\sum_ {j<i}\dfrac{\kappa_{ij}\left(\gamma_{ij}x_{j}^{\ast}(\beta)- 
\gamma_{ji}x_{i}^{\ast}(\beta)\right)^{2}}{ \alpha_{j}\alpha_{i}x_{j}^{\ast}(\beta)x_{i}^{\ast}(\beta)}\geq \sum_{i=1}^{n}K_i.
\end{equation}
Equality  holds if and only if $\beta=0$ or 
$\gamma_{ij}x_{j}^{\ast}(\beta)-\gamma_{ji}x_{i}^{\ast}(\beta)=0$, for all $i$ and $j$. Let us prove that if at least two patches have different growth rates, then equality cannot hold for $\beta>0$. Suppose that there exists $\beta^\ast>0$ such that the positive equilibrium satisfies
\begin{equation}\label{44}
\forall i , j, \quad
\gamma_{ij}x_{j}^{\ast}(\beta^\ast)= \gamma_{ji}x_{i}^{\ast}(\beta^\ast).
\end{equation}
Replacing the equation \eqref{44} in the system \eqref{19}, we get that $x_i^\ast(\beta^\ast)=K_i$, for all $i$. Therefore, from \eqref{44}, it is seen that, for all $i$ and $j$, $K_j\gamma_{ij}=K_i\gamma_{ji}$. From these equations and the conditions 
$\alpha_i\gamma_{ij}=\alpha_j\gamma_{ji}$, we get $r_i=r_j$, for all $i$ and $j$. This is a contradiction with the hypothesis that there exists two patches with different growth rates. Hence the equality in \eqref{som}  holds if and only if $\beta=0$.
\end{proof}

When the matrix $\Gamma$ is irreducible and symmetric, the hypothesis of Prop. \ref{p55} implies that $\alpha_i=\alpha_j$ for all $i$ and $j$. Indeed if two patches $i$ and $j$ are connected (i.e $\gamma_{ij}=\gamma_{ji}\neq 0$), then we have $\alpha_i=\alpha_j$. As the matrix $\Gamma$ is irreducible, for two arbitrary patches, there exists a finite sequence $(i,\ldots, j)$ which begins in $i$ and ends in $j$, such that $\gamma_{ab}\neq 0$ for all successive patches $a$ and $b$ in $(i,\ldots, j)$. Hence  
$\alpha_a=\alpha_b$ for all $a$ and $b$ in $(i,\ldots, j)$. Hence, $\alpha_i=\alpha_j$. So, when the matrix $\Gamma$ is symmetric, Prop. \ref{p55} says that if all $\alpha_i$ are equal, dispersal enhances population growth, which is \cite[Prop. 3.3]{4.1}.

Note that, when $n=2$, Prop \ref{p55}
asserts that if $\alpha_2/\alpha_1=\gamma_{12}/\gamma_{21}$, then $X_T^\ast(\beta)>K_1+K_2$, which is a result of
Arditi et al.  \cite[Prop. 2, item b]{2}. See also Prop. \ref{Prop2patch}, and note that the condition
$\alpha_2/\alpha_1=\gamma_{12}/\gamma_{21}$ implies that $(\gamma_{12},\gamma_{21})\in \mathcal{J}_0$. 

For three patches or more, if the matrix $\Gamma$ does not verify the condition $(\forall i,j, \ \gamma_{ij}=0 \Longleftrightarrow \gamma_{ji}=0)$, then the hypothesis of Prop. \ref{p55}, that for all $j<i$, 
$ \alpha_{i}\gamma_{ij}=\alpha_{j}\gamma_{ji}$ cannot be satisfied. Note that the hypothesis $ \alpha_{i}\gamma_{ij}=\alpha_{j}\gamma_{ji}$ implies that, for all $i=1,\ldots,n$, one has
$$\sum_{j=1}^n\frac{\gamma_{ij}}{\alpha_j}=
\sum_{j=1, j\neq i}^n\frac{\gamma_{ij}}{\alpha_j}-\sum_{j=1, j\neq i}^n\frac{\gamma_{ji}}{\alpha_i}
=\sum_{j=1, j\neq i}^n\frac{\alpha_i\gamma_{ij}-\alpha_j\gamma_{ji}}{\alpha_i\alpha_j}=0.$$
Therefore we can make the following remark:
\begin{rem}
The hypothesis of Prop. \ref{p55} implies that $(\frac{1}{\alpha_1},\ldots, \frac{1}{\alpha_n})^T\in \ker \Gamma$.
\end{rem}

We make the following conjecture:
\begin{con}
If $(\frac{1}{\alpha_1},\ldots, \frac{1}{\alpha_n})^T\in \ker \Gamma$ then 
$$X_T^\ast(\beta)\geq \sum_{i=1}^n K_i,\qquad \text{ for all } \beta\geq 0.$$
\end{con}
This conjecture is true for the particular case of Prop. \ref{p55}. It is also true for two-patch models and for $n$-patch models with symmetric dispersal. It  agrees with  Prop. \ref{3.9}.
\begin{propo}\label{p4.5}
The derivative of the total equilibrium population $X_T^\ast(\beta) $ at $\beta= 0$ is given by:
\begin{equation}\label{deriv}
\frac{dX_T^\ast}{d\beta}(0)=\sum_{i=1}^{n}\left(\frac{1}{r_i}\sum_{j=1}^{n}\gamma_{ij}K_j\right).
\end{equation} 
In particular, if $K\in \ker \Gamma$, where $K=(K_1,\ldots,K_n)^T$, then
$
\frac{dX_T^\ast}{d\beta}(0)=0.
$
\end{propo}
\begin{proof}
By differentiating the equation \eqref{2424} at $\beta=0$, we get:
$$\frac{dX_T^\ast}{d\beta}(0)=\sum_{i=1}^{n}\sum_{j=1 , j \neq i}^{n}\dfrac{\gamma_{ij}x_{j}^{\ast}(0)-\gamma_{ji}x_{i}^{\ast}(0)}{\alpha_{i}x_{i}^{\ast}(0)},$$
which gives \eqref{deriv}, since $x_i^\ast(0)=K_i$ for all $i=1,\ldots,n$.

If $K\in \ker \Gamma$, then $\sum_{j=1 }^{n}\gamma_{ij}K_{j}=0$ for all $i$, so that 
$
\frac{dX_T^\ast}{d\beta}(0)=0.
$
\end{proof}

Actually, when $K\in \ker\Gamma$, we prove that $X_T^\ast(\beta)$ is constant, so that $\frac{dX_T^\ast}{d\beta}(\beta)=0$ for all $\beta\geq0$, not only for $\beta=0$, see Proposition \ref{p52}.
\subsection{Independence of the total equilibrium population with respect to  asymmetric dispersal}
 In the next proposition we give  sufficient and necessary conditions for the total equilibrium population not to depend on the  migration rate.

\begin{propo}\label{p52}
The equilibrium $E^\ast(\beta)$  does not depend on $\beta $  if and only if $(K_{1},\ldots,K_{n})^T\in \ker \Gamma$.
In this case we have $E^\ast(\beta)=(K_{1},\ldots,K_{n})$ for all $\beta >0$.
\end{propo}
\begin{proof}
The equilibrium $E^\ast(\beta)$ is the unique positive solution of the equation 
\begin{equation}\label{4e}
f(x)+\beta\Gamma x=0,
\end{equation}
where $f$ is given by \eqref{f}.
 Suppose that the equilibrium $E^\ast(\beta)$ does not depend on $\beta$, then we replace in Equation \eqref{4e}:
 \begin{equation}\label{4ee}
f(E^\ast(\beta))+\beta\Gamma E^\ast(\beta)=0.
\end{equation}
The derivative of \eqref{4ee} with respect to $\beta$ gives
\begin{equation}\label{4eee}
\Gamma E^\ast(\beta)=0.
\end{equation}
Replacing the equation \eqref{4eee} in the equation \eqref{4ee}, we get $f(E^\ast(\beta))=0$, so $E^\ast(\beta)=(K_{1},\ldots,K_{n})$.
From the equation \eqref{4eee}, we conclude that  $(K_{1},\ldots,K_{n})^T\in \ker \Gamma$.

Now, suppose that $(K_{1},\ldots,K_{n})^T\in \ker \Gamma$, then $(K_{1},\ldots,K_{n})$ satisfies the equation \eqref{4e}, for all $\beta\geq 0$. So, $E^\ast(\beta)= (K_{1},\ldots,K_{n})$, for all $\beta\geq 0$, which proves that the total equilibrium population is independent of the migration rate $\beta$.
\end{proof}

If the matrix $\Gamma$ is symmetric, the previous proposition asserts that the $K_i$, for $i=1,\ldots,n$, are equal if and only if $E^\ast=(K,\ldots, K)$, where $K$ is the common value of the $K_i$. This is \cite[Proposition 3.2]{4.1}. For $n=2$ , Prop. \ref{p52} asserts that if $K_1/K_2=\gamma_{12}/\gamma_{21}$ then $X_T^\ast(\beta)=K_1+K_2$ for all $\beta$, which is \cite[Proposition 2, item c ]{2}. See also the last item of Prop. \ref{Prop2patch}.
\subsection{Two blocks of identical patches}
We consider the model \eqref{m6} and we assume that there are two blocks, denoted $I$ and $J$, of identical patches, such that $I\cup J=\{1,\cdots,n\}$. Let $p$ be the number of patches in $I$ and $q=n-p$ be the number of patches in $J$. 
Without loss of generality we can take
$I=\{1,\cdots,p\}$ and $J=\{p+1,\cdots,n\}$. The patches being identical means that they have the same specific growth rate $r_i$ and carrying capacity $K_i$. Therefore we have
\begin{equation}
\label{samepatches}
\begin{array}{lcl}
r_1=\cdots=r_p,
&&
K_1=\cdots=K_p,\\
r_{p+1}=\cdots=r_n
,&&
K_{p+1}=\cdots=K_n.
\end{array}
\end{equation}

For each patch $i\in I$ we denote by 
$\gamma_{iJ}$ the flux from block $J$ to patch $i$,  and for each patch $j\in J$ we denote by 
$\gamma_{jI}$ the flux from block $I$ to patch $j$, as defined in Table \ref{notations}. For each patch $i$ we denote by $T_i$ the sum of all migration rates 
$\gamma_{ji}$ from patch $i$ to another patch $j\neq i$ (i.e. the outgoing flux of patch i) minus the sum of the migration rates 
$\gamma_{ik}$ from patch $k$ to patch $i$, where $k$ belongs to the same block as $i$. Hence, we have:
\begin{equation}
\label{TiTj}
\left\lbrace \begin{array}{l}
\displaystyle
\mbox{If }i\in I,\quad \mbox{then} \quad
T_i
=\sum_{j\in J}\gamma_{ji}+
\sum_{k\in I\setminus \{i\}}(\gamma_{ki}-\gamma_{ik}).
\\
\displaystyle
\mbox{If }
j\in J,\quad \mbox{then}\quad
T_j
=\sum_{i\in I}\gamma_{ij}+
\sum_{k\in J\setminus \{j\}}(\gamma_{kj}-\gamma_{jk}).
\end{array}\right. 
\end{equation}

We make the following assumption on the migration rates: 
\begin{equation}
\label{samemigrations}
\begin{array}{lcl}
\gamma_{1J}=\cdots=\gamma_{pJ},
&&
\gamma_{(p+1)I}=\cdots=\gamma_{nI}\\
T_{1}=\cdots=T_p
,&&
T_{p+1}=\cdots=T_n
\end{array}
\end{equation}
where $\gamma_{iJ}$, for $i\in I$ and $\gamma_{jI}$, for $j\in J$ are defined in  Table \ref{notations} and $T_i$ are given by (\ref{TiTj}).
\begin{table}[ht]
\caption{Definitions and notations of fluxes\label{notations}}
\begin{tabular}{c|l}
Flux&Definition\\
\hline \\
$
\displaystyle
\gamma_{iJ}=\sum_{j\in J}\gamma_{ij}
$& 
\begin{tabular}{l}
For
$i\in I$,  
$\gamma_{iJ}$ is the flux from block $J$ to patch $i$, i.e. the sum\\ 
of the migration rates $\gamma_{ij}$ from patch $j\in J$ to patch $i$.

\end{tabular}

\\
$\displaystyle
\gamma_{jI}=\sum_{i\in I}\gamma_{ij}
$& 
\begin{tabular}{l}
For
$j\in J$,  
$\gamma_{jI}$ is the flux from block $I$ to patch $j$, i.e. the sum\\ of the migration rates 
$\gamma_{ji}$ from patch $i\in I$ to patch $j$.
\end{tabular}
\\
$\displaystyle
\gamma_{IJ}=
\sum_{i\in I,j\in J}\gamma_{ij}$& 
\begin{tabular}{l}
$\gamma_{IJ}$ is the flux  from block $J$ to block $I$, i.e. the sum\\ of the migration rates 
$\gamma_{ij}$ from patch $j\in J$, to patch $i\in I$.
\end{tabular}
\\
$\displaystyle
\gamma_{JI}=
\sum_{i\in I,j\in J}\gamma_{ji}$& 
\begin{tabular}{l}
$\gamma_{JI}$ is the flux  from block $I$ to block $J$, i.e. the sum\\ of the migration rates 
$\gamma_{ji}$ from patch $i\in I$, to patch $j\in J$.
\end{tabular}
\end{tabular}
\end{table}

We have the following result: 
\begin{lm}
Assume that the conditions (\ref{samemigrations}) are satisfied, then for all $i\in I$ and $j\in I$ one has
\begin{equation}
\label{gammaIJ}
\gamma_{iJ}=\gamma_{IJ}/{p},\quad
\gamma_{jI}=\gamma_{JI}/{q},\quad
T_i=\gamma_{JI}/{p},\quad
T_j=\gamma_{IJ}/{q}.
\end{equation}
where $\gamma_{IJ}$ and $\gamma_{JI}$ are defined in Table \ref{notations}.
\end{lm}
\begin{proof}
The result follows from 
 $\sum_{i\in I}\gamma_{iJ}=\gamma_{IJ}$, 
 $\sum_{i\in J}\gamma_{jI}=\gamma_{JI}$,
 $\sum_{i\in I}T_i=\gamma_{JI}$ and
 $\sum_{i\in J}T_j=\gamma_{IJ}$.
\end{proof}
In the next theorem, we will show that, at the equilibrium, and under certain conditions relating to the  migration rates, we can consider the $n$-patch model as a 2-patch model coupled by migration terms, which are not symmetric in general. Mathematically, we can state our main  result as follows:

\begin{thm}\label{th4}
 Assume that the conditions (\ref{samepatches}) and (\ref{samemigrations}) are satisfied. Then the equilibrium of (\ref{m6}) is of the form 
 $$x_1=x_1^*,\ldots,x_p=x_1^*,
 \quad
 x_{p+1}=x_n^*,\ldots,x_n=x_n^*
 $$ 
where $(x_1^*,x_n^*)$ is the solution of the equations 
\begin{equation}\label{eqx1xnbis}
\left\lbrace \begin{array}{l}
pr_{1}x_{1}\left( 1-\frac{x_{1}}{K_{1}}\right) +
\beta \left({\gamma_{IJ}}x_n-{\gamma_{JI}}x_1\right)=0,
\\
qr_{n}x_{n}\left( 1-\frac{x_{n}}{K_{n}}\right) +\beta \left({\gamma_{JI}}x_1-{\gamma_{IJ}}x_n\right)=0,
\end{array}\right. 
\end{equation}
that is to say, $(x_1^*,x_n^*)$
is the equilibrium  of a 2-patch model, with specific growth rates $pr_1$ and $qr_n$, carrying capacities $K_1$ and $K_n$ and migration rates  
$\gamma_{JI}$ from patch $1$ to patch $2$  and  
$\gamma_{IJ}$ from patch $2$ to patch $1$. 
\end{thm}

\begin{proof}
Assume that the conditions (\ref{samepatches}) are satisfied. Then the equilibrium of (\ref{m6}) is the unique positive solution of the set of algebraic equations
\begin{equation}\label{m6eq}
\left\lbrace \begin{array}{lcl}
\displaystyle
r_{1}x_{i}\left( 1-\frac{x_{i}}{K_{1}}\right) +
\beta \sum_{k=1, k\neq i}^{n}(\gamma_{ik}x_{k}- \gamma_{ki}x_{i})=0,
&& 
i=1,\cdots,p,\\
\displaystyle
r_{n}x_{j}\left( 1-\frac{x_{j}}{K_{n}}\right) +\beta\sum_{k=1, k\neq j}^{n} (\gamma_{jk}x_{k}- \gamma_{kj}x_{j})=0,
&& 
j=p+1,\cdots,n.
\end{array}\right. 
\end{equation}
We consider the following set of algebraic equations obtained from (\ref{m6eq}) by replacing $x_i=x_1$ for $i=1\cdots p$ and 
$x_i=x_n$ for $i=p+1\cdots n$:
\begin{equation}\label{m6eqx1xn}
\left\lbrace \begin{array}{lcl}
r_{1}x_{1}\left( 1-\frac{x_{1}}{K_{1}}\right) +
\beta \left(\gamma_{iJ}x_n-T_ix_1\right)=0,
&& 
i=1,\cdots,p,\\
r_{n}x_{n}\left( 1-\frac{x_{n}}{K_{n}}\right) +\beta \left(\gamma_{jI}x_1-T_jx_n\right)=0,
&& 
j=p+1,\cdots,n.
\end{array}\right. 
\end{equation}
Now, using the assumptions (\ref{samemigrations}), together with the relations (\ref{gammaIJ}), we see that the system (\ref{m6eqx1xn}) is equivalent to the set of two algebraic equations:
\begin{equation}\label{eqx1xn}
\left\lbrace \begin{array}{l}
r_{1}x_{1}\left( 1-\frac{x_{1}}{K_{1}}\right) +
\beta \left(\frac{\gamma_{IJ}}{p}x_n-\frac{\gamma_{JI}}{p}x_1\right)=0,
\\
r_{n}x_{n}\left( 1-\frac{x_{n}}{K_{n}}\right) +\beta \left(\frac{\gamma_{JI}}{q}x_1-\frac{\gamma_{IJ}}{q}x_n\right)=0.
\end{array}\right. 
\end{equation}
We first notice that if 
$x_1 = x_1^*$, $x_n = x_n^*$ is a positive solution of (\ref{eqx1xn}) then
$x_i = x_1^*$ for $i=1,\cdots,p$ and $x_j=x_n^*$ for $j=1,\cdots,n$
is a positive solution of (\ref{m6eq}). Let us prove that (\ref{eqx1xn}) has a unique solution 
$(x_1^*,x_n^*)$. Indeed,  
multiplying the first equation by $p$ and the second one by $q$, we deduce that (\ref{eqx1xn}) can be written in the form 
\eqref{eqx1xnbis}.  
\end{proof}
As a corollary of the previous theorem we obtain the following result which
describes the total equilibrium population in the two blocks:
\begin{cro}\label{Propreduction}
Assume that the conditions (\ref{samepatches}) and (\ref{samemigrations}) are satisfied. Then the total equilibrium population $X_T^*(\beta)=px_1^*(\beta)+qx_n^*(\beta)$ of (\ref{m6}) behaves like the total equilibrium population of the 2-patch model 
\begin{equation}\label{2patchy1yn}
\left\{
 \begin{array}{lcl}
 \frac{dy_1}{dt}&=& r_1y_1\left(1-\frac{y_1}{pK_1}\right) + \beta 
 \left({\gamma_2}{y_n}-{\gamma_{21}}{y_1}\right),\\[6pt]
 \frac{dy_n}{dt}&=& r_ny_n\left(1-\frac{y_n}{qK_n}\right) + \beta
 \left({\gamma_{21}}{y_1}-{\gamma_2}{y_n}\right),
 \end{array} 
\right.
\end{equation}
with specific growth rates $r_1$ and $r_n$, carrying capacities $pK_1$ and $qK_n$, and migration rates   
$\gamma_{21}=\frac{\gamma_{JI}}{p}$, 
$\gamma_2=\frac{\gamma_{IJ}}{q}$.
\end{cro}
\begin{proof}
From Theorem \ref{th4}, we see that
$(x_1^*,x_n^*)$ is the positive solution of (\ref{eqx1xnbis}). 
Hence, $(y_1^*=px_1^*,y_n^*=qx_n^*)$ is the solution of the set of equations
\begin{equation}\label{eqy1yn}
\left\lbrace \begin{array}{l}
r_{1}y_{1}\left( 1-\frac{y_{1}}{pK_{1}}\right) +
\beta \left(\frac{\gamma_{IJ}}{q}y_n-\frac{\gamma_{JI}}{p}y_1\right)=0,
\\
r_{n}y_{n}\left( 1-\frac{y_{n}}{qK_{n}}\right) +\beta \left(\frac{\gamma_{JI}}{p}y_1-\frac{\gamma_{IJ}}{q}y_n\right)=0,
\end{array}\right. 
\end{equation}
obtained from (\ref{eqx1xnbis}) by changing variables to  $y_1=px_1$, $y_n=qx_n$.
The system (\ref{eqy1yn}) has a unique positive solution which is the equilibrium point of the 2-patch model (\ref{2patchy1yn}).  
\end{proof}

We can describe the conditions for which, under Hypothesis (\ref{samepatches}) and (\ref{samemigrations}),  patchiness is beneficial or detrimental in Model (\ref{m6}).

\begin{figure}[ht] 
 \setlength{\unitlength}{0.7cm}
\begin{center}
\begin{picture}(10,7)(0,0.5)
\put(-2,0){\rotatebox{0}{\includegraphics[scale=0.2]{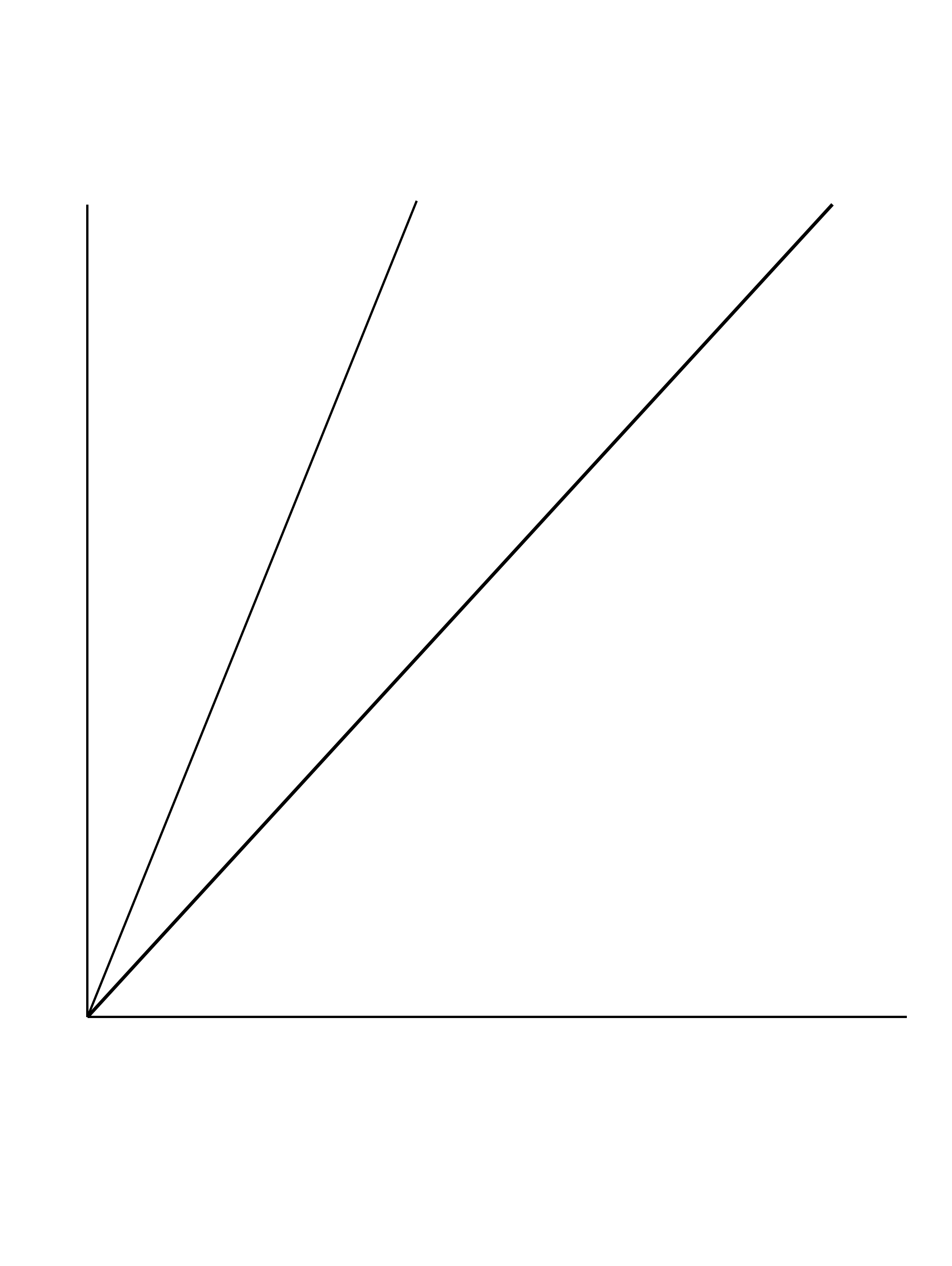}}}
\put(6,0){\rotatebox{0}{\includegraphics[scale=0.2]{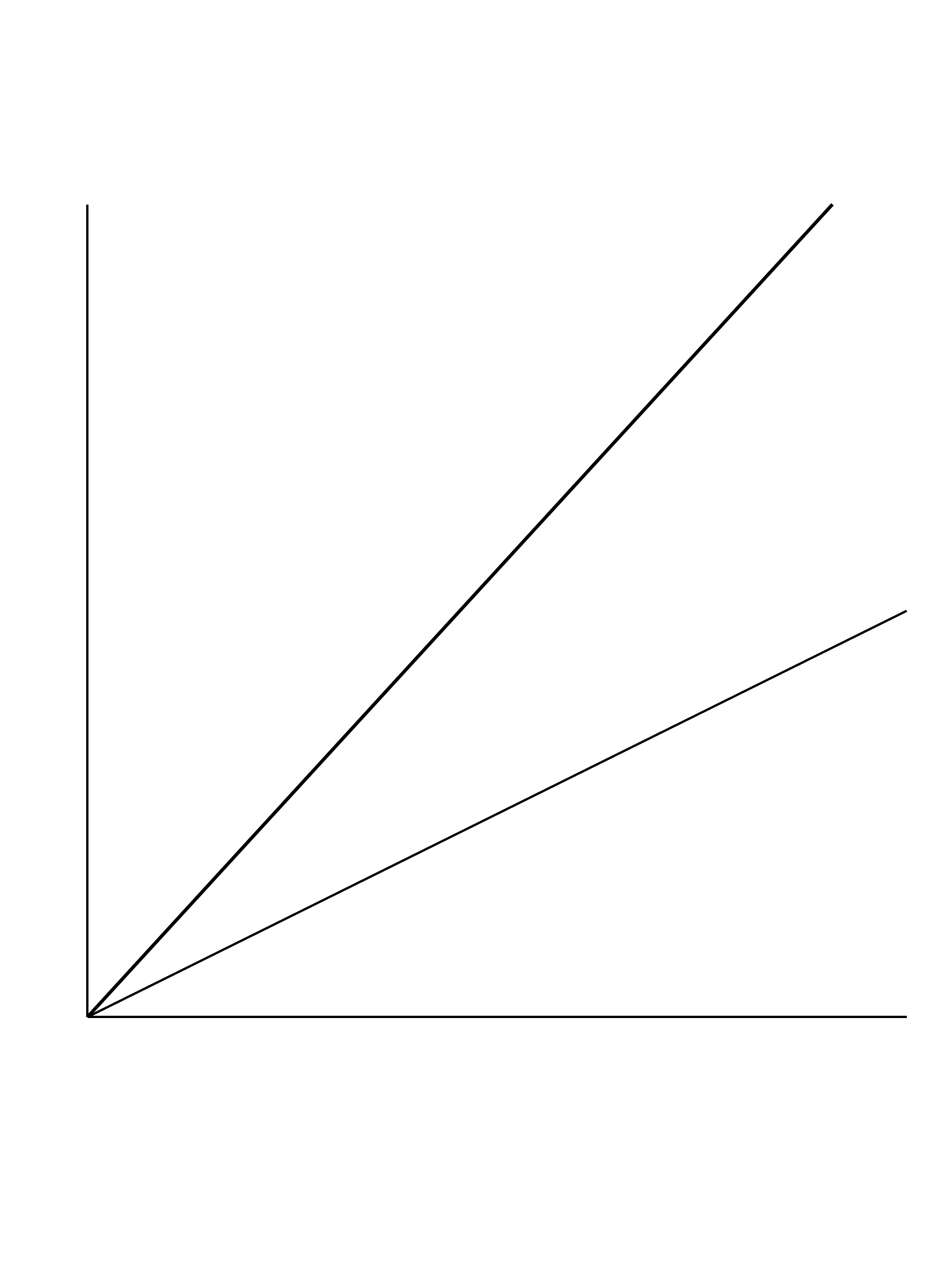}}}
\put(-2,0.5){{Case $r_n>r_1$ (i.e. $\frac{\alpha_n}{\alpha_1}>
\frac{K_1}{K_n}$)}}
\put(-1.8,1.2){$0$}
\put(-1,5){$\mathcal{J}_1$}
\put(0.2,4.4){$\mathcal{J}_0$}
\put(2,3){$\mathcal{J}_2$}
\put(3.8,1.2){$\gamma_{JI}$}
\put(-2,6.7){$\gamma_{IJ}$}
\put(2.4,6.8){{$\frac{\gamma_{IJ}}{\gamma_{JI}}=\frac{K_1}{K_n}$}}
\put(-0.2,6.8){{{$\frac{\gamma_{IJ}}{\gamma_{JI}}
=\frac{\alpha_n}{\alpha_1}$}}}
\put(6,0.5){{Case ${r_n}<{r_1}$
(i.e. $\frac{\alpha_n}{\alpha_1}<
\frac{K_1}{K_n}$)}}
\put(6.2,1.2){$0$}
\put(7.5,5){$\mathcal{J}_2$}
\put(9.3,3.7){$\mathcal{J}_0$}
\put(10,2.5){$\mathcal{J}_1$}
\put(11.8,1.2){$\gamma_{JI}$}
\put(6,6.7){$\gamma_{IJ}$}
\put(10.4,6.8){{$\frac{\gamma_{IJ}}{\gamma_{JI}}=\frac{K_1}{K_n}$}}
\put(10.8,4.2){{{$\frac{\gamma_{IJ}}{\gamma_{JI}}=\frac{\alpha_n}{\alpha_1}$}}}
\end{picture}
\end{center}
   \caption{Qualitative properties of Model (\ref{m6}) under the conditions  (\ref{samepatches}) and (\ref{samemigrations}). In $\mathcal{J}_0$, patchiness has a beneficial effect on the total equilibrium population. This effect is detrimental in $\mathcal{J}_2$. 
In $\mathcal{J}_1$, the effect is beneficial for $\beta<\beta_0$  and detrimental for $\beta>\beta_0$.}
	\label{figurer1r2m6}
\end{figure}
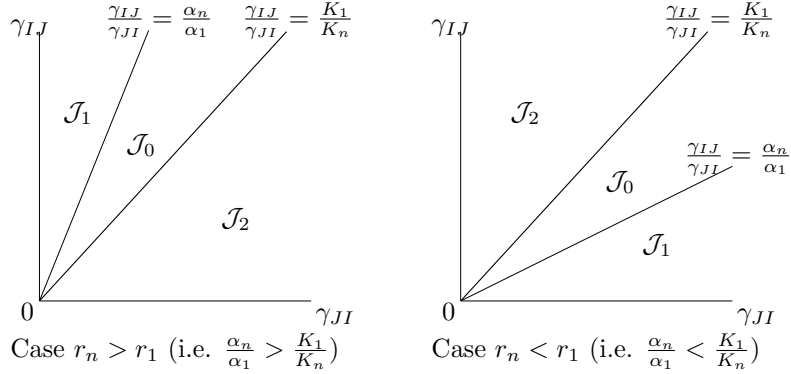
We consider the regions in the set of parameters $\gamma_{IJ}$ and $\gamma_{JI}$, denoted $\mathcal{J}_0$, $\mathcal{J}_1$ and $\mathcal{J}_2$, depicted in 
Fig. \ref{figurer1r2m6} and defined by:
\begin{equation}\label{J0J1J2m6}
\begin{array}{l}
\mbox{If }
r_n>r_1
\mbox{ then }
\left\{
\begin{array}{l}
\mathcal{J}_1=
\left\{(\gamma_{JI},\gamma_{IJ}): \frac{\gamma_{IJ}}{\gamma_{JI}}>\frac{\alpha_n}{\alpha_1}\right\}
\\[6pt]
\mathcal{J}_0=
\left\{(\gamma_{JI},\gamma_{IJ}):\frac{\alpha_n}{\alpha_1}\geq \frac{\gamma_{IJ}}{\gamma_{JI}}>
\frac{K_1}{K_n}\right\}\\[6pt]
\mathcal{J}_2=
\left\{(\gamma_{JI},\gamma_{IJ}): 
\frac{K_1}{K_n}>
\frac{\gamma_{IJ}}{\gamma_{JI}}\right\}
\end{array}
\right.
\\[1cm]
\mbox{If }
r_n<r_1
\mbox{ then }
\left\{
\begin{array}{l}
\mathcal{J}_1=
\left\{(\gamma_{JI},\gamma_{IJ}): \frac{\gamma_{IJ}}{\gamma_{JI}}<\frac{\alpha_n}{\alpha_1}\right\}
\\[6pt]
\mathcal{J}_0=
\left\{(\gamma_{JI},\gamma_{IJ}):\frac{\alpha_n}{\alpha_1}\leq \frac{\gamma_{IJ}}{\gamma_{JI}}<
\frac{K_1}{K_n}\right\}\\[6pt]
\mathcal{J}_2=
\left\{(\gamma_{JI},\gamma_{IJ}): 
\frac{K_1}{K_n}<
\frac{\gamma_{IJ}}{\gamma_{JI}}\right\}
\end{array}
\right.
\end{array}
\end{equation}
where $\alpha_1=r_1/K_1$ and $\alpha_n=r_n/K_n$.

\begin{propo}\label{Propm6}
Assume that the conditions (\ref{samepatches}) and (\ref{samemigrations}) are satisfied. Then the total equilibrium population $X_T^*(\beta)=px_1^*(\beta)+qx_n^*(\beta)$ of  (\ref{m6}) satisfies the following properties
\begin{enumerate}
\item If $r_1=r_n$ then 
$X_T^*(\beta) <  pK_1+qK_n$ for all $\beta > 0$.
\item If ${r_n}\neq{r_1}$, let 
$\mathcal{J}_0$, $\mathcal{J}_1$ and $\mathcal{J}_2$, be defined by (\ref{J0J1J2m6}). Then we have:
\begin{itemize}
\item
if $(\gamma_{JI},\gamma_{IJ})\in \mathcal{J}_0$ then
    $X_T^*(\beta)>pK_1+qK_n$ for any $\beta> 0$,
\item
if $(\gamma_{JI},\gamma_{IJ})\in \mathcal{J}_1$ then
    $X_T^*(\beta)>pK_1+qK_n$ for $0<\beta<\beta_0$  and $X_T^*(\beta)<pK_1+qK_n$ for $\beta>\beta_0$, where
    $$
	\beta_0=\frac{r_n-r_1}
{\frac{\gamma_{IJ}}{\alpha_n}-\frac{\gamma_{JI}}{\alpha_1}}	
	\frac{1}{\frac{\alpha_1}{p}+\frac{\alpha_n}{q}}.$$
\item
if $(\gamma_{JI},\gamma_{IJ})\in \mathcal{J}_2$ then
$X_T^*(\beta)<pK_1+qK_n$ for any $\beta> 0$.
\item If $\frac{\gamma_{IJ}}{\gamma_{JI}}=\frac{K_1}{K_n}$, then 
$X_T^*(\beta)=pK_1+qK_n$ for all $\beta\geq 0$.
 \end{itemize} 
\end{enumerate}
\end{propo}
 
 \begin{proof}
 This is a consequence of Proposition \ref{Prop2patch} and Corollary \ref{Propreduction}.
 \end{proof}
 Let us explain the result of Proposition \ref{Propm6} in the particular case where $p=n-1$. In this case, the condition \eqref{samemigrations} becomes
 \begin{equation}\label{p=n-1}
 \gamma_{1n}=\ldots=\gamma_{n-1,n}\quad  \mbox{ and }\quad T_1=\ldots = T_{n-1},
 \end{equation}
 where $T_i=\gamma_{ni}+\sum\limits_{k\neq i}(\gamma_{ki}- \gamma_{ik})$.\\
 Therefore, if the matrix $\Gamma$ is symmetric, the conditions \eqref{p=n-1} are equivalent to the conditions $ \gamma_{n1}=\ldots=\gamma_{n,n-1}$, which mean that the fluxes of migration between the n-th patch and all $n-1$ identical patches are equal. Hence, Proposition \ref{Propm6}, showing that the n-patch model behaves like a 2-patch model, is the same as  \cite[Prop. 3.4]{4.1}, where the model \eqref{m6} was considered with $\Gamma$ symmetric, $n-1$ patches are identical and the fluxes of migration between the n-th patch and all these $n-1$ identical patches are equal. Thus Proposition \ref{Propm6} generalizes Proposition 3.4 of \cite{4.1}, to asymmetric dispersal and for any two identical blocks, provided that the conditions \eqref{samemigrations} are satisfied.
 \section{Links between SIS and logistic patch models }\label{SECSIS}
 \subsection{The SIS patch model}
In \cite{Gao}, Gao studied the following SIS patch model in an environment of $n$ patches connected by human migration:
\begin{equation}\label{mgao}
\left\lbrace 
\begin{array}{ll}
\frac{dS_i}{dt}=-\beta_i\frac{S_iI_i}{N_i}+\gamma_i I_i+\varepsilon \sum_{j=1}^n \gamma_{ij}S_j, & i=1,\ldots, n,\\[2mm]
\frac{dI_i}{dt}=\beta_i\frac{S_iI_i}{N_i}-\gamma_i I_i+\varepsilon \sum_{j=1}^n \gamma_{ij}I_j, & i=1,\ldots, n,
\end{array}\right. 
\end{equation}
 where $S_i$ and $I_i$ are the number of  susceptible and infected, in patch $i$, respectively; $N_i=S_i+I_i$ denotes the total population in patch $i$. The parameters $\beta_i$ and $\gamma_i$ are positive transmission and recovery rates, respectively. The matrix $\Gamma=(\gamma_{ij})$ satisfies \eqref{gammaiia} and  describes the movement between patches. The coefficient $\varepsilon$ quantifies the diffusion, as our $\beta$ in \eqref{m6}.

 Using the variables $N_i$, $I_i$, $i=1,\ldots,n$, 
 the system \eqref{mgao} has a cascade structure
 \begin{align}\label{mgaoN}
 \frac{dN_i}{dt}&=\varepsilon \sum_{j=1}^n \gamma_{ij} N_i,
 \qquad\qquad\qquad i=1,\ldots, n,\\
 \frac{dI_i}{dt}&=\beta_i\frac{(N_i-I_i)I_i}{N_i}-\gamma_i I_i+\varepsilon \sum_{j=1}^n \gamma_{ij}I_j, \qquad i=1,\ldots, n,
 \end{align}
 Therefore the infected populations $I_i$ are the solutions of the  non-autonomous system of differential equations
\begin{equation}\label{mgaoNA}
 \frac{dI_i}{dt}=\beta_iI_i
 \left(1-\frac{I_i}{N_i(t)}\right)-\gamma_i I_i+\varepsilon \sum_{j=1}^n \gamma_{ij}I_j,\qquad i=1,\ldots, n,
 \end{equation}
 where the total populations $N_i(t)$ are the solutions of the system \eqref{mgaoN}. Hence, the autonomous 
$2n$-dimensional system \eqref{mgao}, is equivalent to the family of $n$-dimensional non-autonomous systems \eqref{mgaoNA}, indexed by the solutions $N_i(t)$ of \eqref{mgaoN}. Note that since the $\gamma_{ij}$ verify the property \eqref{gammaiia}, the total population is constant:
$\sum_{i=1}^n N_i(t)=N$, where 
 $N:=\sum_{i=1}^n \left( S_i(0)+I_i(0)\right)$.
If the matrix $\Gamma=(\gamma_{ij})$ is irreducible, then $N_i(t)$, the total population in patch $i$,  converges towards the limit
\begin{equation}
\label{Ni*}
\lim_{t\to+\infty}N_i(t)=N_i^\ast
\quad\mbox{ where } N_i^\ast:=\frac{N}{\sum_i \delta_i}\delta_i,\qquad i=1,\ldots,n,
\end{equation} 
where $\delta=(\delta_1,\ldots, \delta_n)^T$ is a positive vector which generates the vector space $\ker \Gamma$.
Therefore \eqref{mgaoNA} is an asymptotically autonomous system, whose limit system is obtained by replacing $N_i(t)$ in  \eqref{mgaoNA}, by their limits $N_i^\ast$, given by \eqref{Ni*}:
\begin{equation}\label{gao1}
\frac{dI_i}{dt}=\beta_iI_i\left(1-\frac{I_i}{N_i^\ast} \right)-\gamma_iI_i+\varepsilon \sum_{j=1}^n \gamma_{ij}I_j, \qquad i=1,\ldots, n.
\end{equation} 
The main problem for \eqref{mgao}
is to determine the condition under which the disease free equilibrium, corresponding to the equilibrium $I=0$ of \eqref{gao1}, is GAS, or the endemic equilibrium, corresponding to the positive equilibrium of \eqref{gao1}, is GAS. It is known, see \cite[Theorem 2.1]{Gao}, that the  disease free equilibrium is GAS if $\mathcal{R}_0\leq 1$, and there exists a unique endemic equilibrium, which is GAS, if 
$\mathcal{R}_0> 1$. Here 
$\mathcal{R}_0$ is the basic reproduction number
of the model \eqref{mgao}, defined as:
$$\mathcal{R}_0=\rho\left(FV^{-1} \right)\mbox{ where }
F = {\rm diag}(\beta_1,\cdots,\beta_n)\mbox{ and }V = {\rm diag}(\gamma_1,\cdots,\gamma_n)-\varepsilon \Gamma.$$
A reference work on the basic reproduction number for metapopulations is  Arino \cite{arino}, whereas Castillo-Garsow and  Castillo-Chavez \cite{CGCC} give a
more general account of the subject. 

\subsection{Comparisons between the results on \eqref{m6} and the results on \eqref{gao1}}
Gao \cite{Gao} gave many interesting results on the effect of population dispersal on total
infection size. Our aim is to discuss some of the links between his results and the results of the present paper. We focus on two results on the total infection size $T_n(\varepsilon)=\sum_{i=1}^nI_i^\ast(\varepsilon)$, where $(I_1^\ast(\varepsilon),\ldots,I_n^\ast(\varepsilon))$ is the positive equilibrium of  \eqref{gao1}. We consider the results of Gao \cite{Gao} on 
$T_n(+\infty)$ and 
$T_n'(0)$.
 
\begin{propo}\cite[Theorem 3.3]{Gao}, \cite[Theorem 3.5]{Gao}.
If $\mathcal{R}_0(+\infty)>1$, then
\begin{equation}\label{perfectmixinggao}
T_n(+\infty)=\left( 1-\frac{1}{\mathcal{R}_0(+\infty)}\right) N, \quad \mbox{ with } \quad \mathcal{R}_0(+\infty)=\dfrac{\sum_i \beta_i\delta_i}{\sum_i\gamma_i\delta_i}. 
\end{equation}
If $\beta_i\neq \gamma_i$ for all $i$, then
\begin{equation}\label{dergao}
T_n'(0)=\sum_i\left( \dfrac{1}{\vert\beta_i-\gamma_i\vert}\sum_j \gamma_{ij}I^\ast_j(0)\right),
\quad \mbox{ with } \quad
I^\ast_j(0)=\dfrac{\beta_j-\gamma_j}{\beta_j}N_j^\ast.
\end{equation}
\end{propo}
It is worth noting that the formulas \eqref{perfectmixinggao} and \eqref{dergao} involve the system \eqref{gao1}. An important property of this system is given in the following remark.

\begin{rem}\label{NinkerG}
Let $N^\ast=(N_1^\ast,\ldots, N_n^\ast)^T$ be the vector of the carrying capacities in the system \eqref{gao1}. One has $N^\ast\in \ker\Gamma$, as shown by \eqref{Ni*}.
\end{rem}

Our aim is to compare the results given by the formulas \eqref{perfectmixinggao} and \eqref{dergao} when $\gamma_i\to 0$, to our results, for the system
\begin{equation}\label{gao2}
\frac{dx_i}{dt}=\beta_ix_i\left(1-\frac{x_i}{N_i^\ast} \right)+\varepsilon \sum_{j=1}^n \gamma_{ij}x_j, \qquad i=1,\ldots, n.
\end{equation}
Note that the system \eqref{gao1} reduces to \eqref{gao2} when $\gamma_i= 0$ for all $i$. More precisely we show that, as $\gamma_i\to 0$, the formulas \eqref{perfectmixinggao} and \eqref{dergao} are the same as the results predicted by Prop. \ref{p4.2} and  Prop.\ref{p4.5}.

\begin{propo} Let $T_n(\varepsilon)$ be the total infection size of \eqref{gao1}. Let $X_T^\ast(\varepsilon)$ be the total population size of \eqref{gao2}. One has
\begin{equation}\label{gaoelbetch}
\lim_{\max_i \{\gamma_i\}\to 0}T_n(+\infty)=X_T^\ast(+\infty)=N,\qquad
\lim_{\max_i \{\gamma_i\}\to 0}T_n'(0)=\frac{dX_T^\ast}{d\epsilon}(0)=0.
\end{equation}
\end{propo}\label{p.gaoelbetch}
\begin{proof}
When $\gamma_i\to 0$ for all $i$, one has 
$\mathcal{R}_0(+\infty)\to +\infty$ and 
$I^\ast_i(0)\to N_i^\ast$.
Therefore, from \eqref{perfectmixinggao} and \eqref{dergao} it is deduced that
\begin{equation}\label{resultsgao}
T_n(+\infty)\to N,
\qquad
T_n'(0)\to \sum_i\frac{1}{\beta_i}\sum_j \gamma_{ij}N^\ast_j=0.
\end{equation}
Using the property 
$N^\ast\in\ker\Gamma$, from Prop. \ref{p4.2} and  Prop.\ref{p4.5}, it is deduced that:
\begin{equation}\label{resultselbetch}
X_T^\ast(+\infty)=N,\qquad 
\frac{dX_T^\ast}{d\epsilon}(0)=0.
\end{equation}
From \eqref{resultsgao} and \eqref{resultselbetch} we deduce \eqref{gaoelbetch}.
\end{proof}
Actually as shown in Prop.\ref{p52}, we have the stronger result 
$X_T^\ast(\beta)=N$ for all $\beta\geq 0$. But our aim here was only the comparison between \eqref{resultsgao} and \eqref{resultselbetch}.

As shown in Prop.\ref{p.gaoelbetch}, the results of Gao \cite{Gao} on the logistic patch model \eqref{gao1} yield results on the logistic patch model \eqref{gao2} by taking the limit $\gamma_i\to 0$. However, the scope of this approach is weakened by the fact that it only applies to the logistic model \eqref{gao2}, for which the vector of carrying capacities satisfies $N^\ast\in\ker\Gamma$, see Remark \ref{NinkerG}. But this property is not true in general for our system \eqref{m6}, where the condition $K\in\ker\Gamma$ does not hold in general.

Our aim in this section is to show that any logistic patch model \eqref{m6}, without the condition  $K\in \ker\Gamma$, can be 
written in the form \eqref{gao1}, with the condition 
$N^\ast\in \ker\Gamma$. 
Indeed we have the following result:
\begin{lm}Consider $r_i>0$, $K_i>0$ and $\Gamma$ as in the system \eqref{m6}. Let $\delta_i>0$ be such that $(\delta_1,...,\delta_n)^T\in\ker\Gamma$. Let $N$ be such that $N> \frac{\sum_i \delta_i}{\delta_i}K_i$ for $i=1,\ldots,n$. Let $N_i^\ast$ defined by \eqref{Ni*}. Let $\beta_i=\frac{r_i}{K_i}N_i^\ast$ and $\gamma_i=\beta_i-r_i$. Then one has
\begin{equation}\label{2.1=5.4}
r_ix_i\left(1-{x_i}/{K_i}\right)=
\beta_ix_i\left(1-{x_i}/{N_i^\ast}\right)-\gamma_ix_i,\qquad\mbox{for }i=1,\ldots n
\end{equation}
\end{lm}
\begin{proof}
The conditions \eqref{2.1=5.4} are satisfied if and only if $r_i=\beta_i-\gamma_i$ and 
$r_i/{K_i}=\beta_i/{N_i^\ast}$.
Therefore
\begin{equation}\label{changeparameters}
\left\lbrace 
\begin{array}{l}
\beta_i=N_i^\ast \frac{r_i}{K_i}=N_i^\ast\alpha_i, \\
\gamma_i=\beta_i-r_i=(N_i^\ast-K_i)\alpha_i.
\end{array}
\right. 
\end{equation}
To ensure that $\gamma_i>0$ for all $i$, just choose $N$ in \eqref{Ni*} such that $N_i^\ast>K_i$ for $i=1,\ldots,n$, that is to say, $N>\frac{\sum_i \delta_i}{\delta_i}K_i $.
\end{proof}

\begin{rem}
According to the change of parameters \eqref{changeparameters}, the logistic patch model \eqref{m6} can be written in the form of Gao \eqref{gao1}, i.e. with the property that $N^\ast\in \ker\Gamma$. For the perfect mixing case, the formula \eqref{perfectmixinggao} and our formula \eqref{5151} are the same. Indeed
replacing $\beta_i$ and $\gamma_i$ by \eqref{changeparameters} in \eqref{perfectmixinggao}, and using \eqref{Ni*}, we get:
$$\left( 1-\frac{1}{\mathcal{R}_0(+\infty)}\right) N=\left( 1-\frac{\sum_i (N_i^\ast-K_i)\alpha_i}{\sum_iN_i^\ast\alpha_i} \right)N=\sum_i\delta_i\frac{\sum_i r_i\delta_i}{\sum_i \alpha_i\delta_i^2}.$$
For the derivative, the formula \eqref{dergao} and our formula \eqref{deriv} are the same. Indeed, if we replace $\beta_i$ and $\gamma_i$ by \eqref{changeparameters}, in \eqref{dergao}, we get:
$$I_j^*(0)=\frac{\beta_j-\gamma_j}{\beta_j}N_j^\ast=
\frac{r_j}{N_j^\ast\alpha_j}N_j^\ast=K_j.$$ 
Therefore
$$\sum_i\left( \dfrac{1}{\vert\beta_i-\gamma_i\vert}\sum_j \gamma_{ij}I^\ast_j(0)\right)=\sum_i\left( \dfrac{1}{r_i}\sum_j \gamma_{ij}K_j\right).$$

\end{rem}

The theory of asymptotically autonomous systems  answers  the question ``under which conditions do the solutions of the original $2n$-dimensional system \eqref{mgao} have
the same asymptotic behavior as those of the $n$-dimensional limit system \eqref{gao1} ?''. For details and further reading on the theory of asymptotically autonomous systems the reader is referred to Markus \cite{Markus} and Thieme \cite{Thieme,Thieme1994}. For applications of this theory to epidemic models, see Castillo-Chavez and Thieme \cite{CCT}.

Hence, it is important to know whether or not some of the results of Gao \cite{Gao} on the SIS model  \eqref{mgao} can be deduced from our results on the logistic model \eqref{m6}.
It is worth noting that the discussion in this section shows that our results on the logistic patch model imply results on the model \eqref{gao1} and hence, results on the original model $2n$-dimensional system \eqref{mgao}. However, it is needed that $\beta_i>\gamma_i$ for $i=1,\ldots,n$. Indeed, according to  \eqref{changeparameters}, one has $r_i=\beta_i-\gamma_i>0$. On the other hand, the condition $\beta_i>\gamma_i$ is not required in all patches of the system \eqref{mgao}.  Another challenging problem is the study of the model \eqref{gao1}, in the general case where $N^\ast=(N_1^\ast,\ldots,N_n^\ast)^T$ is not necessarily in the kernel of $\Gamma$.

 \section{Three-patch model}\label{sec5}
\begin{figure}[h]
\begin{center}
\begin{tabular}{cc}
\begin{tikzpicture}
\node[draw,circle](1) at (90:1) {$1$};
\node[draw,circle](2) at (210:1) {$2$};
\node[draw,circle](3) at (-30:1) {$3$};
\draw[>=latex,->] (1) to[bend right] (2);
\draw[>=latex,->] (2) to (1);
\draw[>=latex,->] (1) to(3);
\draw[>=latex,->] (3) to[bend right] (1);
\draw[>=latex,->] (2) to[bend right] (3);
\draw[>=latex,->] (3) to (2);
\end{tikzpicture}
&
\begin{tikzpicture}
\node[draw,circle](1) at (90:1) {$1$};
\node[draw,circle](2) at (210:1) {$2$};
\node[draw,circle](3) at (-30:1) {$3$};
\draw[>=latex,->] (1) to[bend right] (2);
\draw[>=latex,->] (2) to (1);
\draw[>=latex,->] (1) to(3);
\draw[>=latex,->] (3) to[bend right] (1);
\end{tikzpicture}
\\
\\
$\mathcal{G}_1$
&
$\mathcal{G}_2$
\end{tabular}
\end{center}
\caption{The two graphs $\mathcal{G}_1$ and $\mathcal{G}_2$ for which  the  migration matrix may be symmetric, if $\gamma_{ij}=\gamma_{ji}$.
\label{fig1}}
\end{figure}
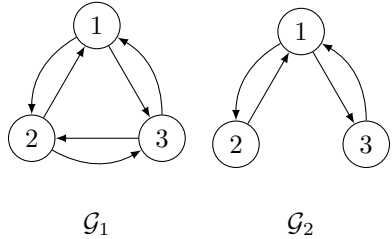

In this section, we consider the model of three patches coupled by asymmetrical terms of migrations. Under the irreducibility hypothesis on the matrix $\Gamma$, there are five possible cases, modulo permutation of the three patches, see Figures \ref{fig1} and \ref{fig2}.

The  connectivity matrices associated to the graphs $\mathcal{G}_1$ and $\mathcal{G}_2$ are given by
\begin{displaymath}
\Gamma_0^{(1)}=\left[\begin{array}{ccc}
0&\gamma_{12}&\gamma_{13}\\
\gamma_{21}&0&\gamma_{23}\\
\gamma_{31}&\gamma_{32}&0\\
\end{array}\right], \quad \text{ and }\quad
\Gamma_0^{(2)}=\left[\begin{array}{ccc}
0&\gamma_{12}&\gamma_{13}\\
\gamma_{21}&0&0\\
\gamma_{31}&0&0\\
\end{array}\right].
\end{displaymath}

For the remaining cases, the graphs $\mathcal{G}_3, \mathcal{G}_4$ and $\mathcal{G}_5$, cannot be symmetrical:
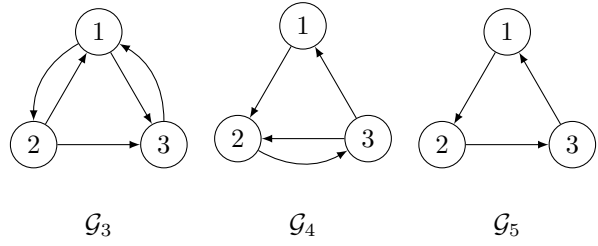
\begin{figure}[h]
\begin{center}
\begin{tabular}{ccc}
\begin{tikzpicture}
\node[draw,circle](1) at (90:1) {$1$};
\node[draw,circle](2) at (210:1) {$2$};
\node[draw,circle](3) at (-30:1) {$3$};
\draw[>=latex,->] (1) to[bend right] (2);
\draw[>=latex,->] (2) to (1);
\draw[>=latex,->] (1) to(3);
\draw[>=latex,->] (3) to[bend right] (1);
\draw[>=latex,->] (2) to (3);
\end{tikzpicture}
&
\begin{tikzpicture}
\node[draw,circle](1) at (90:1) {$1$};
\node[draw,circle](2) at (210:1) {$2$};
\node[draw,circle](3) at (-30:1) {$3$};
\draw[>=latex,->] (1) to (2);
\draw[>=latex,->] (3) to(1);
\draw[>=latex,->] (2) to[bend right] (3);
\draw[>=latex,->] (3) to (2);
\end{tikzpicture}
&
\begin{tikzpicture}
\node[draw,circle](1) at (90:1) {$1$};
\node[draw,circle](2) at (210:1) {$2$};
\node[draw,circle](3) at (-30:1) {$3$};
\draw[>=latex,->] (1) to (2);
\draw[>=latex,->] (3) to(1);
\draw[>=latex,->] (2) to (3);
\end{tikzpicture}
\\
\\
$\mathcal{G}_3$
&
$ \mathcal{G}_4$ 
&
$\mathcal{G}_5$

\end{tabular}
\end{center}
\caption{The three graphs $\mathcal{G}_3, \mathcal{G}_4 $ and $\mathcal{G}_5$ for which the migration matrix cannot be symmetric.\label{fig2}}
\end{figure}

The associated connectivity matrices are given by
\begin{displaymath}
\Gamma_0^{(3)}=\left[\begin{array}{ccc}
0&\gamma_{12}&\gamma_{13}\\
\gamma_{21}&0&0\\
\gamma_{31}&\gamma_{32}&0\\
\end{array}\right],
 \Gamma_0^{(4)}=\left[\begin{array}{ccc}
0&0&\gamma_{13}\\
\gamma_{21}&0&\gamma_{23}\\
0&\gamma_{32}&0\\
\end{array}\right],
\Gamma_0^{(5)}=\left[\begin{array}{ccc}
0&0&\gamma_{13}\\
\gamma_{21}&0&0\\
0&\gamma_{32}&0\\
\end{array}\right].
\end{displaymath}

\begin{figure}[h]
\setlength{\unitlength}{1.0cm}
\begin{center}
\begin{picture}(8.5,5.5)(-2,0.5)		
\put(-4.5,0){{\includegraphics[scale=0.3]{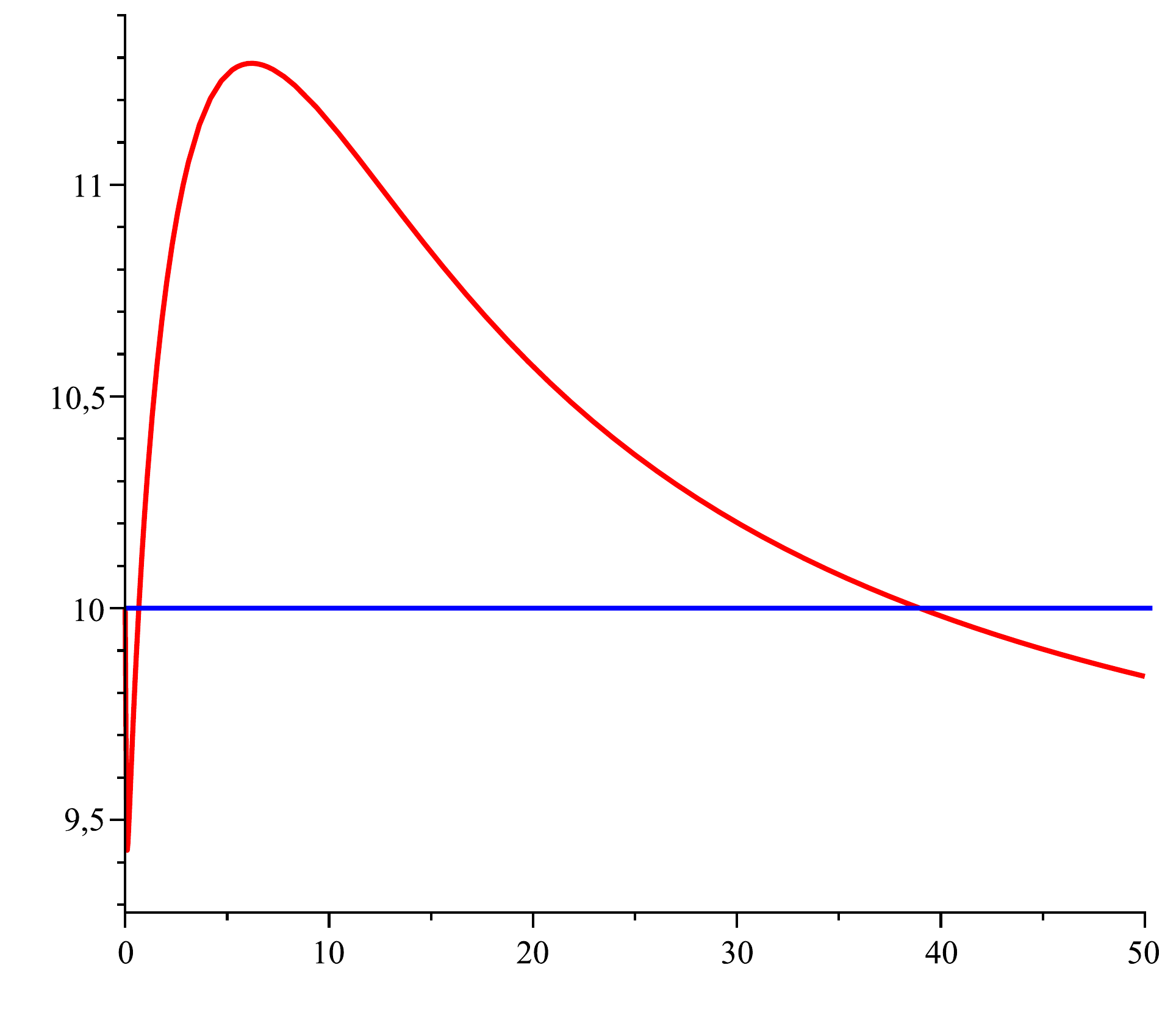}}}
\put(2,0){{\includegraphics[scale=0.29]{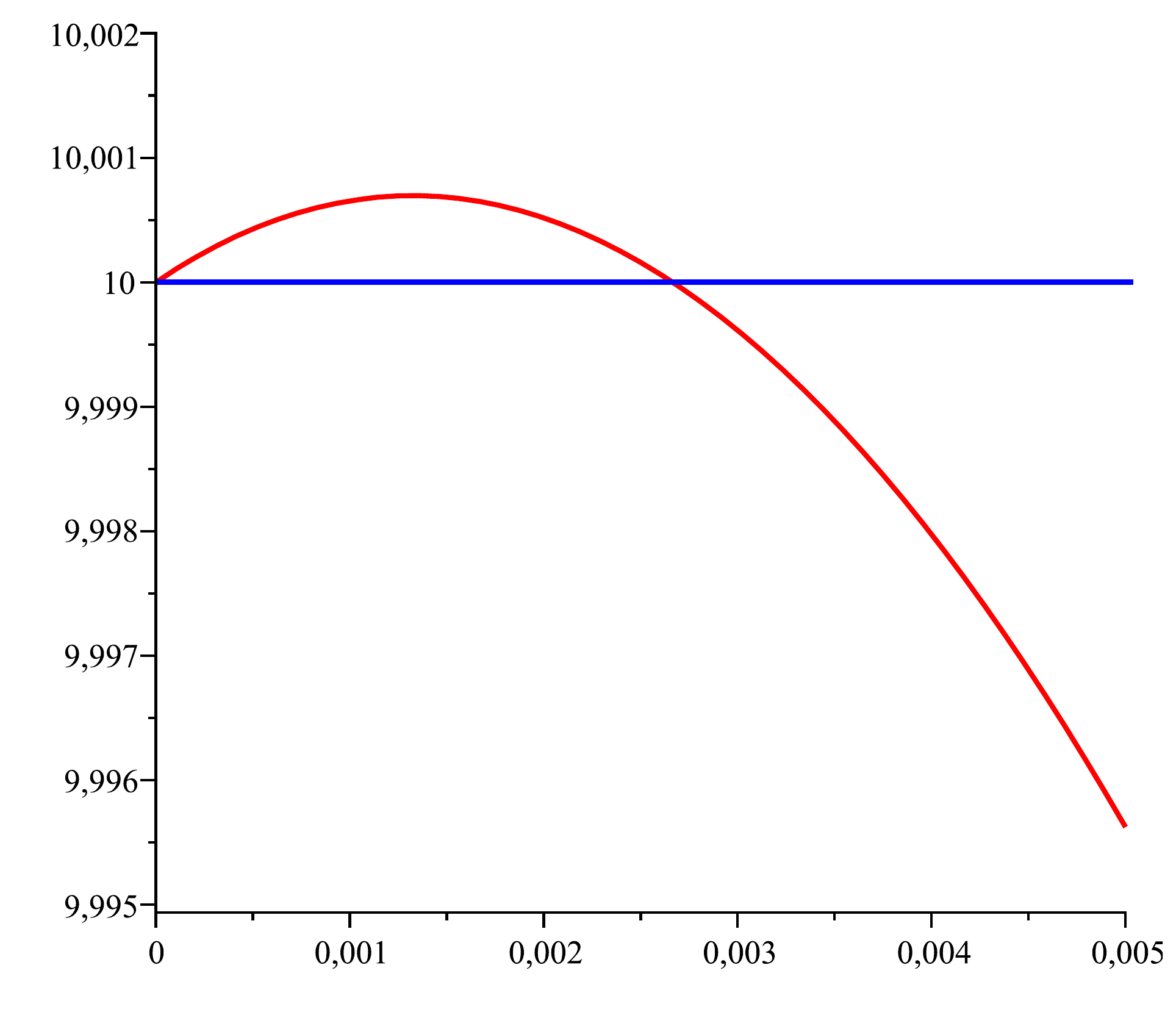}}}
\put(7.7,0.6){$\beta$}
\put(1,0.6){$\beta$}
\put(6.2,3.2){$K_1+K_2+K_3$}
\put(-0.3,2.4){$K_1+K_2+K_3$}		
\put(-4,5.2){$X_T^\ast$}
\put(2.5,5.2){$X_T^\ast$}	
\end{picture}
\end{center}

\caption{Total equilibrium population $X_{T}^{\ast}$ of the system \eqref{m6} $(n=3)$ as a function of  the migration rate $\beta$. The figure on the right is a zoom, near the origin, of the figure on the left. The parameter values are given in Table \ref{t1}.}
\label{f1}
\end{figure}

\begin{figure}[h]
\setlength{\unitlength}{1.0cm}
\begin{center}
\begin{picture}(8.5,6.7)(-2,0.5)		
\put(-4.5,2){{\includegraphics[scale=0.32]{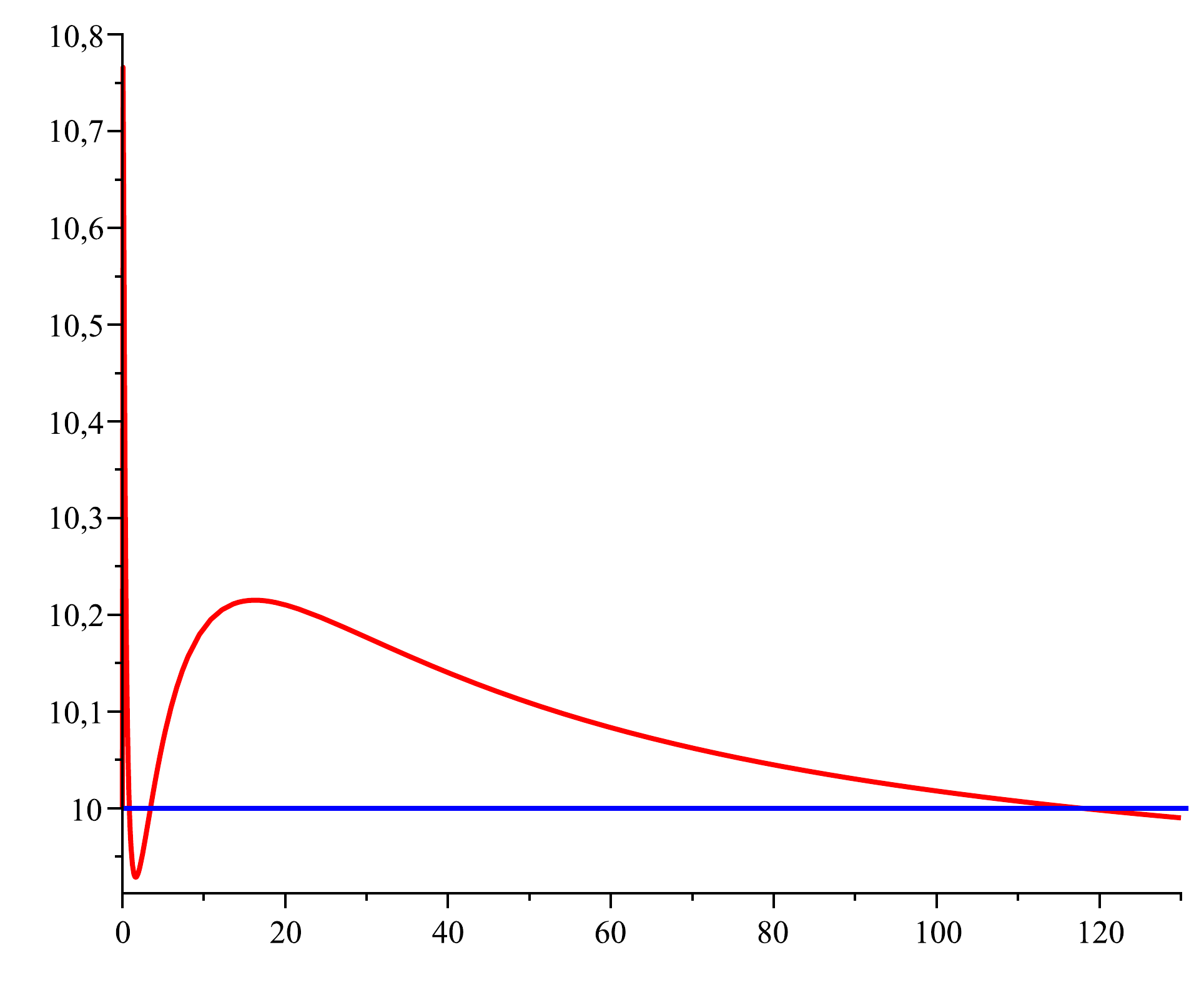}}}                          
\put(2,2){{\includegraphics[scale=0.32]{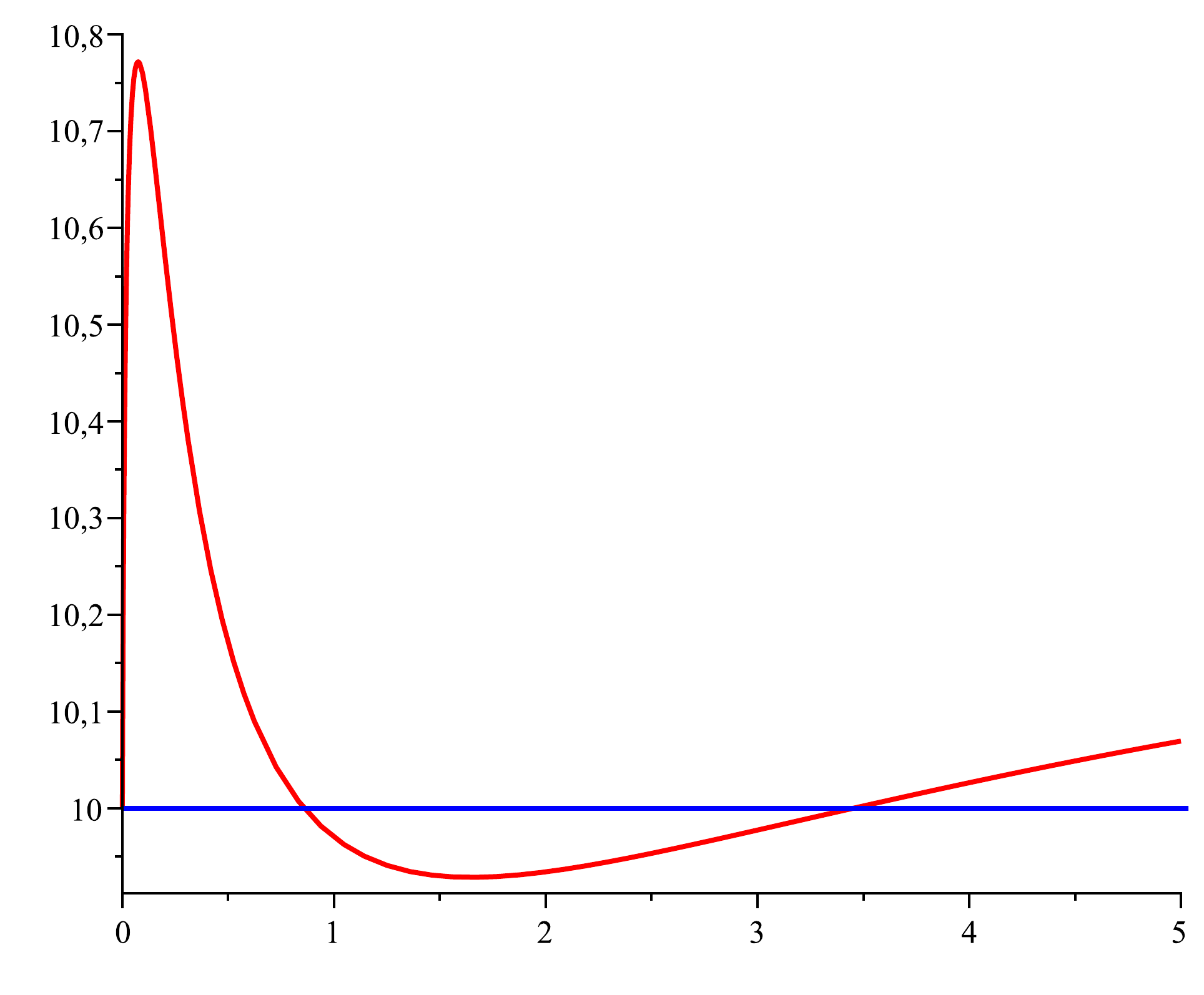}}}                                                 
\put(8,2.1){$\beta$}
\put(1.6,2.1){$\beta$}
\put(-0.3,3.5){$K_1+K_2+K_3$}		
\put(-4,7.1){$X_T^\ast$}
\put(2.5,7.1){$X_T^\ast$}	
\end{picture}
\end{center}

\caption{Total equilibrium population $X_{T}^{\ast}$ of the system \eqref{m6} $(n=3)$ as a function of the migration rate $\beta$. The figure on the right is a zoom, near the origin, of the figure on the left. The parameter values are given in Table \ref{t1}.}
\label{f2}
\end{figure}

\begin{figure}[h]
\setlength{\unitlength}{1.0cm}
\begin{center}
\begin{picture}(8.5,5.5)(-3,0.5)		
\put(-4.5,0){{\includegraphics[scale=0.3]{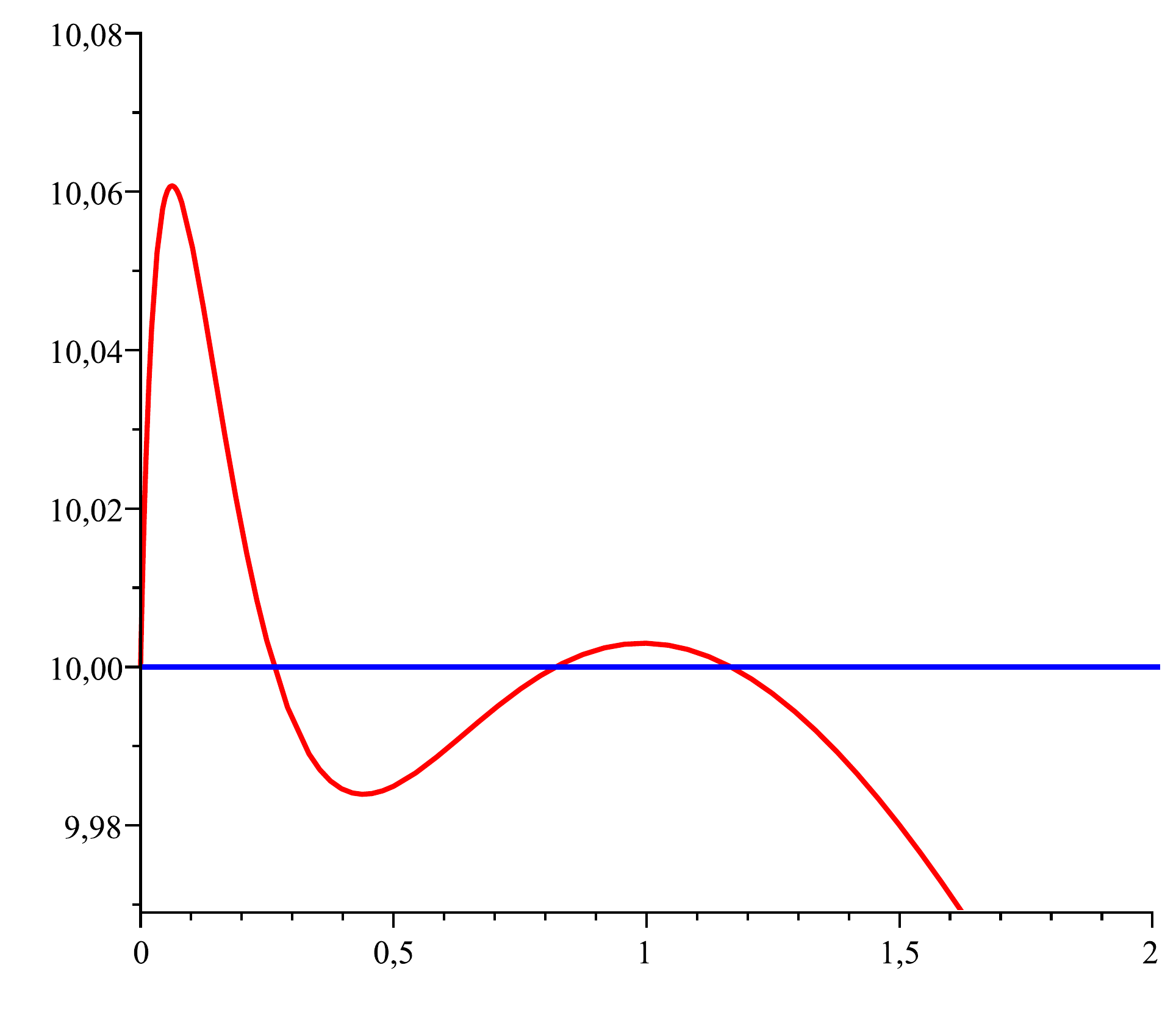}}}
\put(2,0){{\includegraphics[scale=0.3]{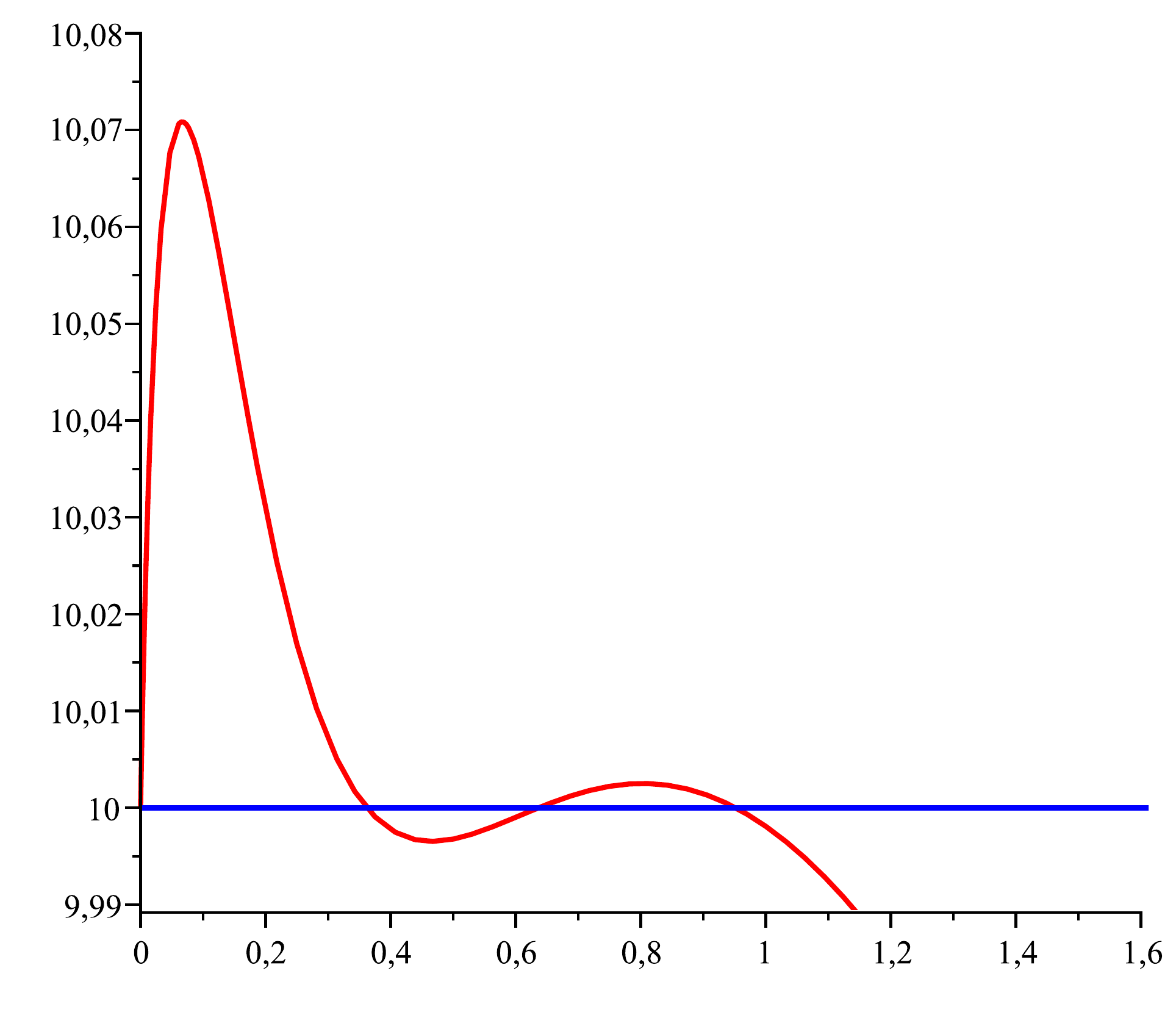}}}
\put(7.7,0.6){$\beta$}
\put(1,0.6){$\beta$}
\put(6.1,1.2){$K_1+K_2+K_3$}
\put(-0.3,1.8){$K_1+K_2+K_3$}		
\put(-4,5.2){$X_T^\ast$}
\put(2.5,5.2){$X_T^\ast$}	
\put(-3,5){$(\textbf{a})$}	
\put(3.5,5){$(\textbf{b})$}	
\end{picture}
\end{center}
\begin{center}
\begin{picture}(8.5,5.5)(-3,0.5)		
\put(-4.5,0){{\includegraphics[scale=0.3]{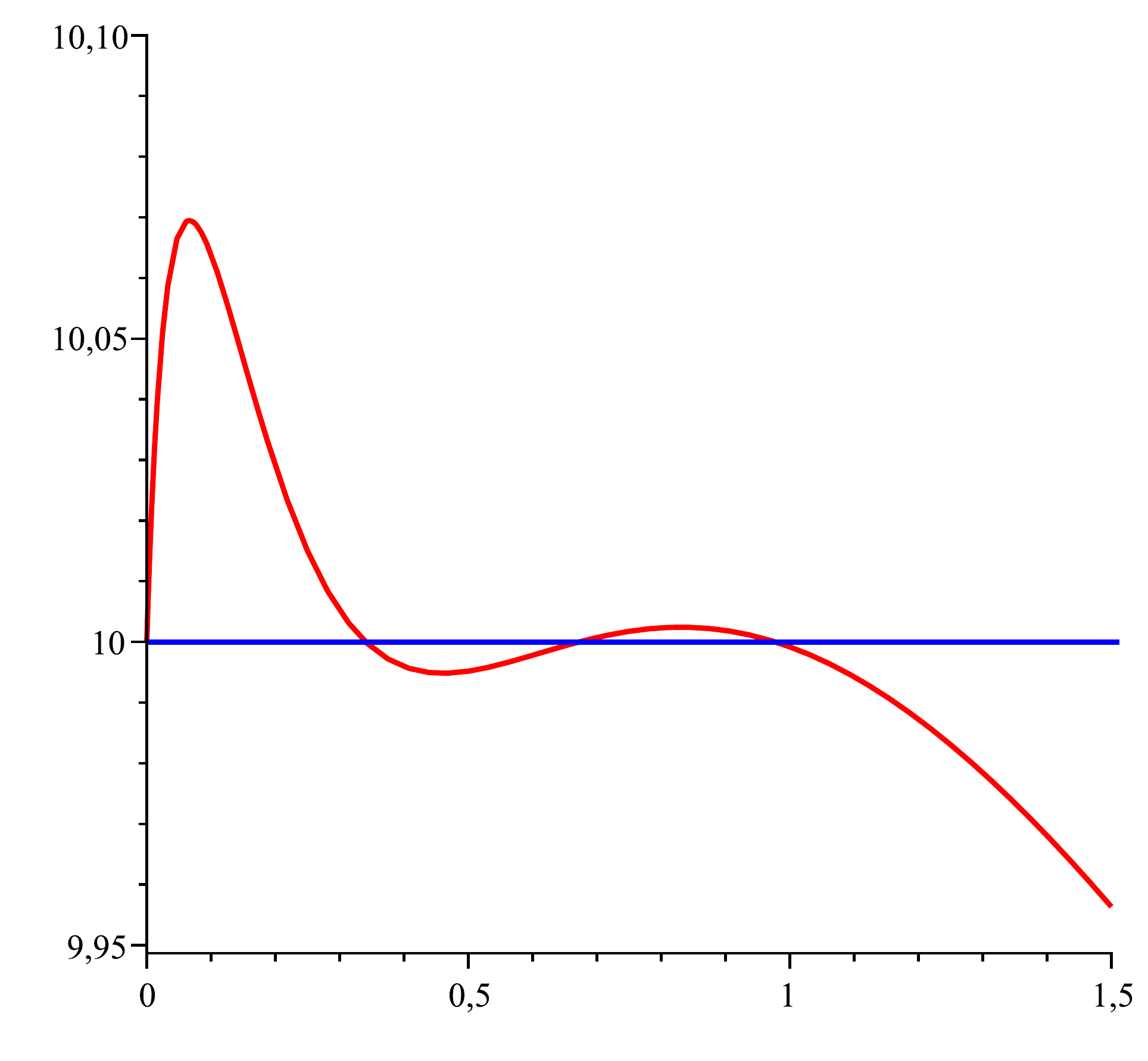}}}

\put(1.2,0.6){$\beta$}
\put(-0.3,2.2){$K_1+K_2+K_3$}		
\put(-4,5.5){$X_T^\ast$}
\put(-3,5){$(\textbf{c})$}	
\end{picture}
\end{center}

\caption{Total equilibrium population $X_{T}^{\ast}$ of the system \eqref{m6} $(n=3)$ as a function of the migration rate $\beta$. 
The parameter values are given in Table \ref{t1}.\label{f3}}
\end{figure}
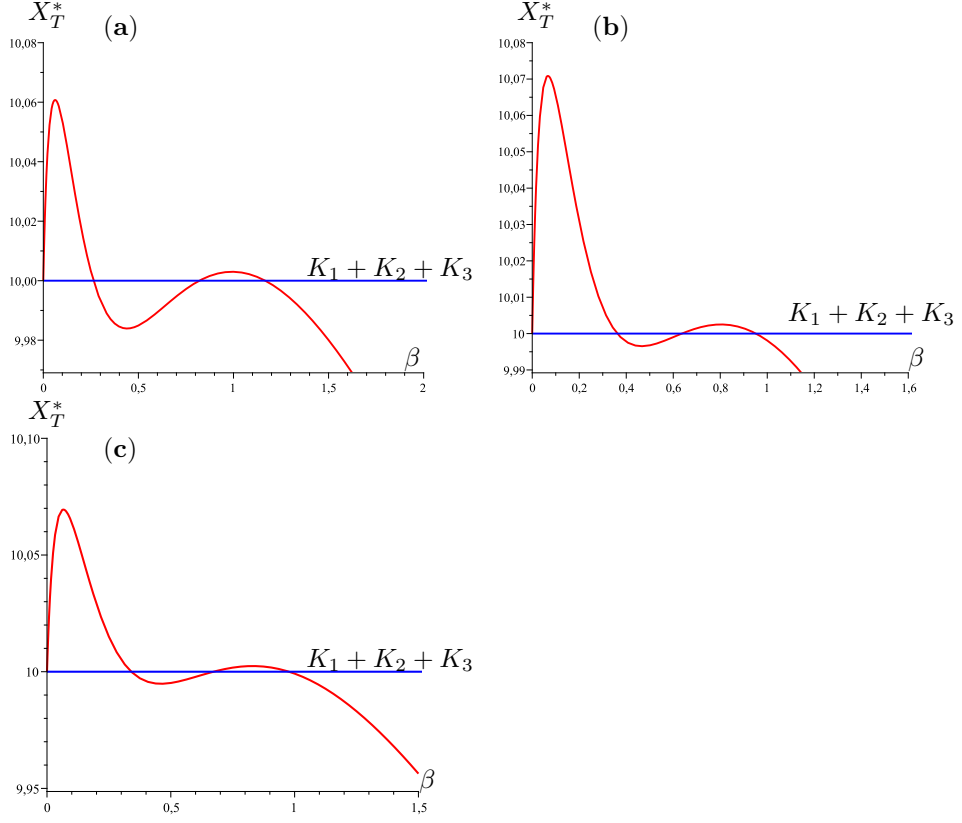

In Table \ref{tt1}, we give the formula of perfect mixing $X_T^\ast(+\infty)$ for each of the  five cases.
\begin{small}
\begin{table}[h]
\caption{ The generator $\delta$ of $\ker\Gamma$, for the five cases. The perfect mixing abundance $X_{T}^{\ast}(+\infty)$ is computed with Eq. \eqref{5151}.\label{tt1}}
\begin{tabular}{l|l}
\hline 
\\
Graphs &The formula of perfect mixing $X_{T}^{\ast}(+\infty)$ \\ 
\hline \\
$\mathcal{G}_1$& The coefficients $\delta_i$ are given by the equation \eqref{coedel}\\ \\
$\mathcal{G}_2$& $\delta_1=\gamma_{12}\gamma_{13}, \delta_2=\gamma_{{21}}\gamma_{{13}}, \delta_3=\gamma_{{12}}\gamma_{{31}},$   \\ \\
$\mathcal{G}_3$  & $\delta_1=\gamma_{12}\gamma_{13}+\gamma_{32}\gamma_{13}, 
\delta_2=\gamma_{21}\gamma_{13}, 
\delta_3=\gamma_{21}\gamma_{32}+\gamma_{31}\gamma_{12}+\gamma_{31}\gamma_{32},$\\ \\
$\mathcal{G}_4$  &$\delta_1=\gamma_{32}\gamma_{13},\delta_2=\gamma_{21}\gamma_{13}+\gamma_{21}\gamma_{23}+\gamma_{31}\gamma_{23},\delta_3=\gamma_{21}\gamma_{32}.$\\ \\ 
$\mathcal{G}_5$& $\delta_1=\gamma_{32}\gamma_{13},\delta_2=\gamma_{21}\gamma_{13},
\delta_3=\gamma_{21}\gamma_{32}.$\\ \\  
\hline 
\end{tabular} 
\end{table}
\end{small}
In the numerical simulations, we show that we can have new behaviors of $X_T^* (\beta)$. 
In the case $n=2$, it was shown in \cite{1,2} that there exists at most one positive value of $\beta$ such that  $X_T^* (\beta)=K_1+K_2$.
 In \cite{4.1},  in the case $n=3$ and $\Gamma$ is symmetric, we gave numerical values for the parameters such that there exists two positive values of $\beta$ such that  
 $X_T^* (\beta)=K_1+K_2+K_3$, and we were not able to find more than two values. 
 The novelty when $\Gamma$ is not symmetric is that we can find examples with three positive values. 
 Indeed, we may have the following situation :
 $\frac{dX_{T}^{\ast}}{d\beta}(0)>0$ and 
$X_{T}^{\ast}(+\infty)< K_{1}+K_{2}+K_{3}$, and there exist three values  $0<\beta_1<\beta_2<\beta_3$ for which we have\begin{equation}\label{e2}
X_{T}^{\ast}(\beta)    \left\lbrace 
\begin{array}{lll}
>K_{1}+K_{2}+K_{3}& \text{ for }& \beta\in \left]0,\beta_1\right[    \cup \left]\beta_2,\beta_3 \right[,\\
< K_{1}+K_{2}+K_{3}& \text{ for }& \beta\in  \left]\beta_1,\beta_2 \right[  \cup  \left]\beta_3,+\infty \right[.
\end{array} \right.
\end{equation}
 The same situation  holds for each of the five  graphs $\mathcal{G}_1$, $\mathcal{G}_2, \mathcal{G}_3,\mathcal{G}_4$ and $\mathcal{G}_5$, i.e, there exist three values  $0<\beta_1<\beta_2<\beta_3$ for which \eqref{e2} hold. See Figures \ref{f1},  (for the graph $\mathcal{G}_1$), \ref{f2}, (for the graph $\mathcal{G}_2$), 
 \ref{f3}-a, (for the graph $\mathcal{G}_3$),
 \ref{f3}-b, (for the graph $\mathcal{G}_4$),
 and \ref{f3}-c, (for the graph $\mathcal{G}_5$).

\begin{table}[h]
\caption{
The numerical values of the parameters for the logistic growth function and migration coefficients of the model \eqref{m6}, with $n=3$, used in Fig.~\ref{f1},\ref{f2},\ref{f3}-a,\ref{f3}-b and Fig \ref{f3}-c. For all  figures  we have $(r_1, r_2, r_3, K_1, K_2, K_3)=(4, 0.7, 0.6, 5, 1, 4)$. 
The perfect mixing abundance $X_{T}^{\ast}(+\infty)$ is computed with Eq. \eqref{5151} and the derivative of the total equilibrium population at $\beta= 0$ is computed with Eq. \eqref{deriv}\label{t1}.}
\begin{tabular}{l|cccccccccccccc}
\hline \\
Figure & $\gamma_{21}$&$\gamma_{12}$&$\gamma_{31}$&$\gamma_{13}$&$\gamma_{32}$&$\gamma_{23}$& $\frac{dX^{\ast}_{T}}{d \beta}(0)$&$X_{T}^{\ast}(+\infty)$ \\ 
\hline

 \ref{f1}&$0.15$&$3$&$0.2$ &  $0.04$&$11$&$0.1$&  $1.06$& $9.21$\\ 
 \ref{f2}&$14.9$&$10$&$0.2$ &  $0.04$&$\textcolor{red}{0}$&$\textcolor{red}{0}$& $77.20$& $9.86$\\ 
 \ref{f3}-a&$1.44$&$0.01$&$0.2$ &  $0.04$&$1$&$\textcolor{red}{0}$ & $3.11$& $8.93$\\ 
\ref{f3}-b&$1.52$&$\textcolor{red}{0}$&  $\textcolor{red}{0}$&$1$&$1$ &$0.002$& $3.52$& $8.72$\\ 
\ref{f3}-c&$1.51$&$\textcolor{red}{0}$&  $\textcolor{red}{0}$&$1$&$1$ &$\textcolor{red}{0}$ & $3.46$& $8.75$\\  
\end{tabular} 
\end{table}

\section{Conclusion}
The aim of this paper is to generalize, to a multi-patch model with asymmetric dispersal, the results obtained in \cite{4.1} for a multi-patch model with symmetric dispersal. 
 
In Section \ref{sec3} we consider the particular case of perfect mixing, when the migration rate goes to infinity, that is, individuals may travel freely between patches. 
As in \cite{4.1}, we compute the total equilibrium population in that case, and, by perturbation arguments, 
we prove  that the dynamics in this ideal case provides a good approximation to the case when the migration rate is large. 
Our results generalize those of \cite{2} (asymmetric migration matrix, only two patches), \cite{3.2} 
(arbitrarily many patches, but the migration matrix is symmetric and zero outside the corners and the  three main diagonals), 
and \cite{4.1} (arbitrarily many patches;  arbitrary, but symmetric, migration).
 
In Section \ref{sec4} we consider the equation 
\begin{equation}\label{E}
\mbox{ total equilibrium population }=\mbox{ sum of the carrying capacities of the patches}.
\end{equation} 
We give a complete solution in the case when the $n $ patches are partitioned into two blocks of identical patches. 
Our results mirror those of \cite{2}, which deals with the two-patch case. Specifically, Equation (\ref{E}) admits at most one non-trivial solution. 
 
 In Section \ref{SECSIS}, we consider a SIS patch model and we give the links with the logistic model.
 
In Section \ref{sec5} we give numerical values for the dispersion parameters such that Equation (\ref{E})   has at least three non-trivial solutions. 
In \cite{4.1} we proved that for three patches and symmetric dispersal, there may be at least two solutions. 
A mathematical proof that, when n=3, Equation (\ref{E})  has at most three solutions,  would certainly be desirable, and could spur further work. 
Upper bounds for arbitrarily many patches would also be interesting. 
\section*{Appendix}
\appendix
\begin{normalsize}
 \section{The 2-patch asymmetric model}\label{apa}
We consider the 2-patch logistic equation with asymmetric migrations. 
We denote by $\gamma_{12}$ the migration rate from patch 2 to patch 1 and $\gamma_{21}$ from patch 1 to patch 2. The model is written: 

\begin{equation}\label{2patch}
\left\{
 \begin{array}{lcl}
 \dfrac{dx_1}{dt}&=& r_1x_1\left(1-\dfrac{x_1}{L_1}\right) + \beta 
 \left({\gamma_{12}}{x_2}-{\gamma_{21}}{x_1}\right),\\[6pt]
 \dfrac{dx_2}{dt}&=& r_2x_2\left(1-\dfrac{x_2}{L_2}\right) + \beta
 \left({\gamma_{21}}{x_1}-{\gamma_{12}}{x_2}\right).
 \end{array} 
\right.
\end{equation}
Note that the system \eqref{2patch} is  studied in \cite{1,3,5,6,9} in  the case where the migration rates satisfy  $\gamma_{21}=\gamma_{12}$,  and in \cite{2} for  general migration rates. 
This system admits a unique equilibrium which is GAS. We denote by $E^\ast(\beta)=(x_1^\ast(\beta), x_2^\ast(\beta))$ this equilibrium and by 
$X_T^\ast(\beta)$ the sum of $x_i^\ast(\beta)$.
\begin{figure}[h]
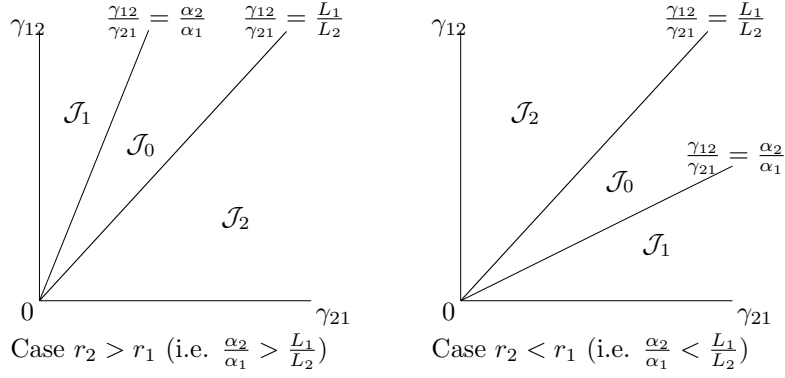
 
 \setlength{\unitlength}{0.7cm}
\begin{center}
\begin{picture}(10,7)(0,0.5)
\put(-2,0){\rotatebox{0}{\includegraphics[scale=0.2]{Fgamma1gamma21.pdf}}}
\put(6,0){\rotatebox{0}{\includegraphics[scale=0.2]{Fgamma1gamma22.pdf}}}
\put(-2,0.5){{Case $r_2>r_1$ (i.e. $\frac{\alpha_2}{\alpha_1}>
\frac{L_1}{L_2}$)}}
\put(-1.8,1.2){$0$}
\put(-1,5){$\mathcal{J}_1$}
\put(0.2,4.4){$\mathcal{J}_0$}
\put(2,3){$\mathcal{J}_2$}
\put(3.8,1.2){$\gamma_{21}$}
\put(-2,6.7){$\gamma_{12}$}
\put(2.4,6.8){{$\frac{\gamma_{12}}{\gamma_{21}}=\frac{L_1}{L_2}$}}
\put(-0.2,6.8){{{$\frac{\gamma_{12}}{\gamma_{21}}=\frac{\alpha_2}{\alpha_1}$}}}
\put(6,0.5){{Case ${r_2}<{r_1}$
(i.e. $\frac{\alpha_2}{\alpha_1}<
\frac{L_1}{L_2}$)}}
\put(6.2,1.2){$0$}
\put(7.5,5){$\mathcal{J}_2$}
\put(9.3,3.7){$\mathcal{J}_0$}
\put(10,2.5){$\mathcal{J}_1$}
\put(11.8,1.2){$\gamma_{21}$}
\put(6,6.7){$\gamma_{12}$}
\put(10.4,6.8){{$\frac{\gamma_{12}}{\gamma_{21}}=\frac{L_1}{L_2}$}}
\put(10.8,4.2){{{$\frac{\gamma_{12}}{\gamma_{21}}=\frac{\alpha_2}{\alpha_1}$}}}
\end{picture}
\end{center}
   \caption{Qualitative properties of model (\ref{2patch}). In $\mathcal{J}_0$, patchiness has a beneficial effect on total equilibrium population. This effect is detrimental in $\mathcal{J}_2$. 
In $\mathcal{J}_1$, the effect is beneficial for $\beta<\beta_0$  and detrimental for $\beta>\beta_0$. In the figure $\alpha_1=r_1/L_1$ and $\alpha_2=r_2/L_2$.}
	\label{figurer1r2}
\end{figure}
We consider the regions in the set of the parameters $\gamma_{21}$ and $\gamma_{12}$, denoted $\mathcal{J}_0$, $\mathcal{J}_1$ and $\mathcal{J}_2$, depicted in 
Fig. \ref{figurer1r2} and defined by:
\begin{equation}\label{J0J1J2}
\left\{
\begin{array}{l}
\mbox{If }
r_2>r_1
\mbox{ then }
\left\{
\begin{array}{l}
\mathcal{J}_1=
\left\{(\gamma_{21},\gamma_{12}): \frac{\gamma_{12}}{\gamma_{21}}>\frac{\alpha_2}{\alpha_1}\right\}
\\[6pt]
\mathcal{J}_0=
\left\{(\gamma_{21},\gamma_{12}):\frac{\alpha_2}{\alpha_1}\geq \frac{\gamma_{12}}{\gamma_{21}}>
\frac{L_1}{L_2}\right\}\\[6pt]
\mathcal{J}_2=
\left\{(\gamma_{21},\gamma_{12}): 
\frac{L_1}{L_2}>
\frac{\gamma_{12}}{\gamma_{21}}\right\}
\end{array}
\right.
\\[1cm]
\mbox{If }
r_2<r_1
\mbox{ then }
\left\{
\begin{array}{l}
\mathcal{J}_1=
\left\{(\gamma_{21},\gamma_{12}): \frac{\gamma_{12}}{\gamma_{21}}<\frac{\alpha_2}{\alpha_1}\right\}
\\[6pt]
\mathcal{J}_0=
\left\{(\gamma_{21},\gamma_{12}):\frac{\alpha_2}{\alpha_1}\leq \frac{\gamma_{12}}{\gamma_{21}}<
\frac{L_1}{L_2}\right\}\\[6pt]
\mathcal{J}_2=
\left\{(\gamma_{21},\gamma_{12}): 
\frac{L_1}{L_2}<
\frac{\gamma_{12}}{\gamma_{21}}\right\}
\end{array}
\right.
\end{array}
\right.
\end{equation}

We have the following result which gives the conditions for which patchiness is beneficial or detrimental in model (\ref{2patch}).

\begin{propo}\label{Prop2patch}
The total equilibrium population of (\ref{2patch}) satisfies the following properties
\begin{enumerate}
\item If $r_1=r_2$ then 
$X_T^*(\beta)\leq L_1+L_2$ for all $\beta\geq 0$.
\item If ${r_2}\neq{r_1}$, let 
$\mathcal{J}_0$, $\mathcal{J}_1$ and $\mathcal{J}_2$, be defined by (\ref{J0J1J2}). Then we have:
\begin{itemize}
\item
if $(\gamma_{21},\gamma_{12})\in \mathcal{J}_0$ then
    $X_T^*(\beta)>L_1+L_2$ for any $\beta> 0$
\item
if $(\gamma_{21},\gamma_{12})\in \mathcal{J}_1$ then
    $X_T^*(\beta)>L_1+L_2$ for $0<\beta<\beta_0$  and $X_T^*(\beta)<L_1+L_2$ for $\beta>\beta_0$, where
    $$
	\beta_0=\frac{r_2-r_1}
{\frac{\gamma_{12}}{\alpha_2}-\frac{\gamma_{21}}{\alpha_1}}	
	\frac{1}{\alpha_1+\alpha_2}.$$
\item
if $(\gamma_{21},\gamma_{12})\in \mathcal{J}_2$ then
$X_T^*(\beta)<L_1+L_2$ for any $\beta> 0$
\item If $\frac{\gamma_{12}}{\gamma_{21}}=\frac{L_1}{L_2}$, then $x_1^*(\beta)=L_1$ and $x_2^*(\beta)=L_2$ for all $\beta\geq 0$.\\
Therefore $X_T^*(\beta)=L_1+L_2$ for all $\beta\geq 0$.
 \end{itemize} 
\end{enumerate}
\end{propo}
\begin{proof}
This result was established by Arditi et al. \cite{2}. Part (1) is Proposition 1 of \cite{2}. The first three items of part (2) are Proposition 2 of \cite{2}. For the last item of part (2), see the last paragraph in page 12 of \cite{2}. The explicit expression of $\beta_0$ was not given in \cite{2}, however, it is easy to deduce it from the formulas given in \cite{2}. 
\end{proof}


\section{Some useful results}\label{apb}

We begin with a 
\begin{lm}\label{appBlm1}
The matrix $\mathcal{L}$ defined by \eqref{0410} 
is stable, that is to say, all its eigenvalues have negative real part.
\end{lm}
\begin{proof}
We consider the two matrices 
$$
G:=\left[ 
\begin{array}{cccc}
&L-U&& V\\
0&\ldots &0 &0
\end{array}
\right],  \qquad 
P:=\left[ 
\begin{array}{cccc}
&I&& 0\\
1&\ldots &1 &1
\end{array}
\right], 
$$
where  $L$, $V$, and $U$  are defined right after \eqref{0410}.
We prove that the two matrices $\Gamma$ and $G$ are conjugate by the matrix $P$, that is to say $P^{-1}GP= \Gamma.$\\
The inverse of matrix $P$ is given by
$$
P^{-1}=\left[ 
\begin{array}{cccc}
&I&& 0\\
-1&\ldots &-1 &1
\end{array}
\right]. 
$$
We have
$$
P^{-1}G P=\left[ 
\begin{array}{cccc}
&L&& V\\
\gamma_{n1}&\ldots &\gamma_{n n-1} & -\sum_{j=1, j\neq 1}^{n}\gamma_{jn}
\end{array}
\right] =\Gamma.
$$
Two conjugate matrices have the same eigenvalues. As the matrix $G$ is block-triangular, its eigenvalues are zero and the eigenvalues of the matrix $L-U$. 
Therefore,  since $0$ is an simple eigenvalue of the matrix $\Gamma$,  the eigenvalues of the matrix $L-U$ are the eigenvalues of the matrix $\Gamma$ except $0$. 
By Lemma \ref{lm41}   all non-zero eigenvalues of $\Gamma$ have negative real part.
\end{proof}

\begin{lm}\label{lma41}
Let $(u_n)_{n\geq 1}$, $(v_n)_{n\geq 1}$ and $(w_n)_{n\geq 1}$ be three real and non-negative sequences.  Then,
\begin{enumerate}
\item if $(u_n)_{n\geq 1}$ and $(v_n)_{n\geq 1}$ are both  non-increasing, or both non-decreasing, then we have, for all  $N\geq1$,
\begin{equation}\label{a17}
\left(\sum_{n=1}^N w_n \right) \left(\sum_{n=1}^Nw_n u_n v_n\right)\geq \left(\sum_{n=1}^Nw_n u_n\right)\left(\sum_{n=1}^Nw_n v_n \right),
\end{equation}
\item if $(u_n)_{n\geq 1}$ is non-decreasing and $(v_n)_{n\geq 1}$ is non-increasing, or if $(u_n)_{n\geq 1}$ is non-increasing and $(v_n)_{n\geq 1}$ is non-decreasing,  then,  we have, for all  $N\geq1$,
\begin{equation}\label{a18}
\left(\sum_{n=1}^Nw_n \right) \left(\sum_{n=1}^Nw_n u_n v_n\right)\leq \left(\sum_{n=1}^Nw_n u_n\right)\left(\sum_{n=1}^Nw_n v_n \right).
\end{equation}
\end{enumerate}
In both items, if $(u_n)_{n\geq 1}$ is not constant, then the inequality in the conclusion is strict. 
\end{lm}
\begin{proof}
We prove Item 1 by induction on $N$, in the case when $(u_n)_{n\geq 1}$ and $(v_n)_{n\geq 1}$ are both  non-decreasing, the other case being identical. Obviously, Equation \eqref{a17} holds for $N = 1$. Now,
assume that \eqref{a17} holds for $N$, then we proceed to show that
\eqref{a17} holds for $N+1$. Since
$$
u_{n+1}\left[ w_1(v_{n+1}-v_1)+\ldots+w_n (v_{n+1}-v_n) \right] \geq u_1w_1 (v_{n+1}-v_1)+\ldots+ u_n w_n (v_{n+1}-v_n),
$$
the inequality being strict if $(u_n)_{n\geq 1}$ is not constant, we observe that
\begin{equation}\label{a19}
\sum_{n=1}^Nw_n u_n v_n+\left(\sum_{n=1}^Nw_n \right)u_{N+1}v_{N+1}\geq \left(\sum_{n=1}^Nw_n v_n \right)u_{N+1}+\left(\sum_{n=1}^Nw_n u_n \right)v_{N+1}.
\end{equation}

From the induction hypothesis and the equation \eqref{a19}, it follows that
\begin{align*}
\left(\sum_{n=1}^{N+1}w_n \right) \left(\sum_{n=1}^{N+1}w_n u_n v_n\right)&=\left(\sum_{n=1}^Nw_n \right) \left(\sum_{n=1}^Nw_n u_n v_n\right)+w_{N+1} \left(\sum_{n=1}^Nw_n u_n v_n\right)\\
&+w_{N+1}^{2}u_{N+1}v_{N+1}+\left(\sum_{n=1}^{N}w_n \right) w_{N+1}u_{N+1}v_{N+1}\\
&\geq \left(\sum_{n=1}^Nw_n \right) \left(\sum_{n=1}^Nw_n u_n v_n\right)+w_{N+1}^{2}u_{N+1}v_{N+1}\\
&+\left(\sum_{n=1}^Nw_n v_n \right)u_{N+1}w_{N+1}+\left(\sum_{n=1}^Nw_n u_n \right)v_{N+1}w_{N+1}\\
&\geq \left(\sum_{n=1}^Nw_n u_n\right)\left(\sum_{n=1}^Nw_n v_n \right)+w_{N+1}^{2}u_{N+1}v_{N+1}\\
&+\left(\sum_{n=1}^Nw_n v_n \right)u_{N+1}w_{N+1}+\left(\sum_{n=1}^Nw_n u_n \right)v_{N+1}w_{N+1}\\
&=  \left(\sum_{n=1}^{N+1}w_n u_n\right)\left(\sum_{n=1}^{N+1}w_n v_n \right).
\end{align*}
This completes the proof of item 1.\\
Equation \eqref{a18}  can then  be proved by reversing all the inequalities in the proof of \eqref{a17} above.
\end{proof}
This result is proved by DeAngelis et al. \cite[Lemma 2.6]{3.1} for Part (2) and in \cite[Proposition A.3]{3.2} for part (1), where $w_n=1$ for all $n\geq 1$. Here we generalize this result to any positive sequence.

\end{normalsize}

\end{document}